\documentclass[12pt,openbib]{amsart}

\usepackage[english]{babel}
\usepackage[utf8]{inputenc}

\usepackage[a4paper,%
 top = 2cm,%
 bottom = 2cm,%
 left = 2cm,%
 right = 2cm,%
 marginparwidth = 1.75cm%
 ]{geometry}

\usepackage{latexsym}
\usepackage{amsmath}
\usepackage{amsthm}
\usepackage{amsfonts}
\usepackage{amssymb}
\usepackage{mathtools}
\usepackage{esint}
\usepackage{verbatim}
 \allowdisplaybreaks
\usepackage{pstricks}

\usepackage{bbm}

\usepackage{xcolor}

\usepackage{hyperref}

\usepackage{lipsum}

\usepackage{mathtools}
\newcommand{\be}{\begin{equation}}
\newcommand{\ee}{\end{equation}}
\newcommand{\bes}{\begin{equation}}
\newcommand{\ees}{\end{equation}}
 
\mathtoolsset{showonlyrefs=true}

\def\lf{\left}
\def\rg{\right}

\def\al{\alpha}
\def\la{\lambda}

\def\om{\omega}
\def\p{\partial}

\def\11{1\!\!1}
 
\newtheorem{Th}{Theorem}[section]
\newtheorem{Prop}[Th]{Proposition} 
\newtheorem{Co}[Th]{Corollary} 
\newtheorem{Lm}[Th]{Lemma} 
 
\newtheorem{Rm}[Th]{Remark} 
 
\newtheorem{Th*}{Theorem}
\newtheorem{Co*}[Th*]{Corollary}

\def \R{\mathbb{R}}

\def \N{\mathbb{N}}
\def \Z{\mathbb{Z}}
 
\title{Morse Index Stability for the Ginzburg-Landau Approximation}%

\author{Francesca Da Lio, Matilde Gianocca}
\address{Department of Mathematics, ETH Zentrum,
CH-8093 Z\"urich, Switzerland.}
\email{francesca.dalio@math.ethz.ch, matilde.gianocca@math.ethz.ch}
\begin{document}

\begin{abstract}
\noindent In this paper we study the behaviour of critical points of the Ginzburg-Landau perturbation of the Dirichlet energy into the sphere
\begin{equation}
E_\varepsilon(u):=\int_\Sigma \frac{1}{2}|du|^2_h\ \,dvol_h +\frac{1}{4\varepsilon^2}(1-|u|^2)^2\,dvol_h=\int_{\Sigma}e_{\varepsilon}(u)\,.
\end{equation}
Our first main result is a precise point-wise estimate for $e_\varepsilon(u_k)$ in the regions where compactness fails, which also implies the $L^{2,1}$ quantization in the bubbling process. Our second main result consists in applying the method developed in the joint paper \cite{DLGR} with T. Rivière to study the upper-semi-continuity of the extended Morse index to sequences of critical points of $E_{\epsilon}$: given a sequence of critical points $u_{\varepsilon_k}:\Sigma\to \mathbb{R}^{n+1}$ of $E_\varepsilon$ that converges in the bubble tree sense to a harmonic map $u_\infty\in W^{1,2}(\Sigma,{S}^{n})$ and bubbles $v^i_{\infty}:\mathbb{R}^2\to {S}^{n}$, we show that the extended Morse indices of the maps $v^i,u_\infty$ asymptotically control the extended Morse index of the sequence $u_{\varepsilon_k}$ for $k$ large enough.
\end{abstract}

\maketitle

\setcounter{equation}{0}
\setcounter{page}{1}
\setcounter{section}{1}
\section*{Introduction}
In the paper \cite{DLGR}, the authors in collaboration with T. Rivière have developed a new method to show the stability of the extended Morse index --the
number of negative or null eigenvalues-- of sequences of critical points of conformally invariant Lagrangians of the form:
\be
\label{ConfEnergy}
\mathfrak{L}(u):=\frac{1}{2}\int_\Sigma |du|^2_h\ \,dvol_h+u^\ast\al
\ee
where $ u\colon(\Sigma,h)\to({\mathcal{N}}^n,g)$, $(\Sigma,h)$ is a smooth closed and oriented Riemannian surface, $({\mathcal{N}}^n,g)$ is a $C^2$ closed Riemannian manifold of arbitrary dimension $n$, and $\al$ is an arbitrary $2$-form on ${\mathcal{N}}^n$. 
The main result of the joint work \cite{DLGR} is
\begin{Th*}
\label{th-morse-stability}
Let $\al$ be a $C^2$ $2$-differential form of ${\mathcal{N}}^n$ and 
$u_k$ be a sequence of critical points of \eqref{ConfEnergy} from a closed oriented surface $(\Sigma,h)$ into a closed arbitrary at least $C^2$ Riemannian manifold ${\mathcal{N}}^n$ such that {\bf bubble tree converges} towards $(u_\infty, v_\infty^1\ldots v_\infty^Q)$. Then, for $k$ large enough, one has
\begin{equation}
\label{0.5}
\mbox{\rm Ind}_{\mathfrak{L}}(u_k)+\mbox{\rm Null}_{\mathfrak{L}}(u_k)\le \mbox{ \rm Ind}_{\mathfrak{L}}(u_\infty)+\mbox{\rm Null}_{\mathfrak{L}}(u_\infty)+\sum_{j=1}^Q\mbox{\rm Ind}_{\mathfrak{L}}(v_\infty^j)+\mbox{\rm Null}_{\mathfrak{L}}(v_\infty^j). ~~~~ 
\end{equation}
\end{Th*}
We recall that the {\it Morse index} and the {\it nullity} are respectively the dimension of the largest space on which the second derivative of $\mathfrak{L}$ is negative definite and the dimension of the kernel of the bi-linear form corresponding to $D^2\mathfrak{L}$.
Classical results (see \cite{MoRe} or more recently \cite{KaSt} in the case $\alpha=0$) imply that for $k$ large enough one has the lower-semi-continuity of the Morse index in the following sense
\begin{equation}
\label{0.6}
\mbox{Ind}_{\mathfrak{L}}(u_\infty)+\sum_{j=1}^Q\mbox{Ind}_{\mathfrak{L}}(v_\infty^j)\le \mbox{Ind}_{\mathfrak{L}}(u_k)\ .
\end{equation}
 
The lower semi-continuity of the Morse index is a rather general and ``robust'' inequality. It holds in various frameworks like for instance in the viscosity method for minimal surfaces, \cite{Riv2}.

The upper-semi-continuity of the extended Morse index has been established in two remarkable works on minimal surface theory respectively by O. Chodosh and C. Mantoulidis for sequences of critical points to the Allen-Cahn functional (see \cite{ChMa}, Theorem 1.9)
and by F.C. Marques and A. Neves for limits of Almgren-Pitts minmax procedures (see \cite{MaNe}, Theorem 1.3). In both cases a main assumption is that the limiting minimal surface has multiplicity one. It has been disproved when the limit has multiplicity higher than one (see \cite{dPKWY} and Example 5.2 in \cite{ChMa} ). 
Compared to the lower semi-continuity Morse index estimates, upper semi-continuity results are
typically much more complicated due to the need of obtaining precise estimates of the sequence of
solutions in regions of loss of compactness and in fact it does not always hold.
The techniques that we have introduced in \cite{DLGR} to prove the upper-semi-continuity of the extended Morse index associated to conformally invariant variational problems in $2$-D have turned to be very efficient in several other geometrical analysis settings (see \cite{Mi, MiRi2, MiRi3, Work}).
\par
 
Let ${S}^{n}$ denote the $n$-dimensional sphere. We recall that a map $u\in W^{1,2}(\Sigma, {S}^{n})$ is said to be \textit{weakly harmonic} if it is a critical point of the Dirichlet energy with respect to smooth variations in the target ${S}^{n}$, i.e. for 
\begin{equation}\label{Denergy}
 {E}(u)=\frac{1}{2}\int_{\Sigma}|du|^2_h\,dvol_h\,,
\end{equation}
it must hold
\begin{equation}
 \frac{d}{dt} {E}\bigg(\frac{u+t\varphi}{|u+t\varphi|}\bigg)\bigg\vert_{t=0}=0\quad\text{for all}\quad\varphi\in C^{\infty}(\Sigma)\,.
\end{equation}
The Euler-Lagrange equation associated to such critical points is 
\begin{equation}\label{EulerEq}
 - \Delta_h u = |d u|_h^2\,u \quad\text{in}\quad \mathcal{D}'(\Sigma),
\end{equation}
where $\Delta_h$ is denoting the negative Laplace Beltrami operator.\footnote{We recall that the negative Laplace Beltrami operator associated to the metric $h$ acting on a smooth function $f\colon\Sigma\to \mathbb{R}$, is in local coordinates given by
 $$ \Delta_h f=\frac{1}{\sqrt{{\rm det}(h_{ij})}}\partial_{x_i}(h^{ij}{\rm det}(h_{ij})\partial_{x_j}f)$$
 where we use the Einstein summation convention. We also have in local coordinates $|df|_h=h^{ij}\partial_{x_i}f\partial_{x_j}f$ }
The equation is critical in the sense that for solutions $u\in W^{1,2}(\Sigma, {S}^{n})$ the right hand side $|d u|_h^2\,u$ is only in $L^1$ and one cannot directly apply Calderon-Zygmund theory to get a better regularity of the solutions $u$. 
To prove the continuity of harmonic maps taking values into the sphere, F. H\'elein \cite{Hel} identified a Jacobian structure in the Euler-Lagrange equation \eqref{EulerEq} and was able to locally express the Euler-Lagrange equation in the form
\begin{equation}\label{harmsp2}
-\Delta u=\nabla^\perp B\cdot \nabla u,~~~\mbox{in $\mathbb{D}^2$}
\end{equation}
where for simplicity we consider the unit ball $\mathbb{D}^2$ centered at the origin and $\nabla^\perp B=(\nabla^\perp B_{ij})$ with $\nabla^\perp B_{ij}=u_i\nabla u_j-u_j\nabla u_i,$ (for every
vector field $v\colon \mathbb{R}^2\to \mathbb{R}^m$, $\nabla^\perp v$ denotes the $\pi/2$ rotation of the gradient $\nabla v$, namely $\nabla^\perp v=(-\partial_y v,\partial_x v))\,.$\par
The r.h.s of \eqref{harmsp2} can be written actually as a sum of Jacobians:
$$
\nabla^\perp B_{ij}\nabla u_j=\partial_xu_j\partial_y B_{ij}-\partial_yu_j\partial_x B_{ij}\,.
$$
This particular structure allows one to apply Wente's Inequality (\cite{Wen}) to equation \eqref{EulerEq}:
\begin{Th*}
\label{th-I.1}
 Let $a$ and $b$ be two measurable functions in $ \mathbb{D}^2$ whose gradients are in $L^2( \mathbb{D}^2)$. Let Then $\varphi\in W_0^{1,1}( \mathbb{D}^2)$ be the unique solution to
\begin{equation}\label{jac}
\left\{\begin{array}{ll}
- \Delta\varphi =\displaystyle{\frac{\partial a}{\partial x} \frac{\partial b}{\partial y} -\frac{\partial a}{\partial y} \frac{\partial b}{\partial x}},&~\mbox{in $ \mathbb{D}^2$}\\ \noindent
\varphi=0 &~\mbox{on $\partial \mathbb{D}^2\,.$}
\end{array}\right.
\end{equation}
Then $u\in W_0^{2,1}( \mathbb{D}^2)$ and there is a constant $C>0$ independent of $a$ and $b$ such that
$$
||\varphi||_\infty+||\nabla\varphi||_{L^{2,1}}\le C||\nabla a||_{L^2}||\nabla b||_{L^2}\,.
$$
In particular $\varphi$ is a continuous in $ \mathbb{D}^2\,.$~~$\Box$
\end{Th*}
We observe that the continuity easily follows thanks to the embedding $W^{1,(2,1)}\xhookrightarrow{} C^0$ in two dimensions.\\ \noindent

One of the most interesting questions related to the study of the energy \eqref{ConfEnergy} is the existence of non-trivial critical points. Unfortunately the energy \eqref{ConfEnergy} does not satisfy the Palais-Smale condition and therefore one cannot directly apply the classical methods in the calculus of variations. 
 In order to overcome this difficulty one considers suitable relaxations of the energy satisfying the Palais-Smale condition.
 We mention that J. Sacks and K. Uhlenbeck in the ground-breaking paper \cite{SaUh} considered the following regularization of the Dirichlet Energy
\begin{equation}\label{Ep}
E_p(u):=\int_{\Sigma}(1+|\nabla u|^2)^{\frac{p}{2}}\ dx^2
\end{equation}
 where $p>2.$ Since the energy \eqref{Ep} satisfies the Palais-Smale condition, they could produce a sequence of smooth critical points to \eqref{Ep}, called $p$-harmonic maps, which {\bf bubble tree converges} towards $(u_\infty, v_\infty^1\ldots v_\infty^Q)$ and as an application of their blow-up analysis they were able to show the existence of a minimal
 two-sphere under suitable topological conditions on the target manifold. T. Lamm \cite{Lam} later proved that there is no energy loss during the blow-up process under an ad-hoc entropy condition. In a joint recent paper \cite{DLRp} of the first author with T. Rivi\`ere, the authors proved an $L^{2,1}$ estimate result in neck-regions for a sequence of critical points of the energy 
\begin{equation}
\label{Ep2}
E_p(u):=\int_{\Sigma}(1+|\nabla u|^2)^{\frac{p}{2}}\ dx^2 +\int_{\Sigma} u^*\alpha
\end{equation}
under an entropy condition which is a bit stronger than the one used by T. Lamm in \cite{Lam}, and which links the exponent $p$ and the concentration radii appearing in the neck region analysis. The energy identity and the "necklessness" property for $p$-harmonic maps into spheres have been obtained in \cite{LiZh} by J. Li and X. Zhu, by combining the existence of conservation laws with the use of Lorentz spaces, using an approach similar to the one introduced twenty years earlier by F.-H. Lin and T. Rivi\`ere in \cite{LiRi1} and further developed by P. Laurain and T. Rivi\`ere in \cite{LaRi}.
Contrary to the case of Ginzburg-Landau type perturbations that we will present below, energy can be lost and necks can arise during the bubbling process for $p$-harmonic maps when the target manifold is not a sphere. An example of this phenomenon has been found by Y. Li and Y. Wang in \cite{LiWa}.\\ \noindent
 
In this paper we consider a perturbation of \eqref{ConfEnergy}
 given by the Ginzburg-Landau energy
\begin{equation}\label{Eeps}
 {E}_\varepsilon(u)= \int_{\Sigma}\frac{1}{2}|d u|_h^2+\frac{1}{4\varepsilon^2}(1-|u|^2)^2\,dvol_h
\,,
\end{equation}
\par

The energy $ {E}_\varepsilon$ is defined on the Hilbert space $W^{1,2}(\Sigma,\mathbb{R}^{n})$ and can be used to approximate harmonic maps $u:\Sigma\to {S}^{n}$. The analysis of the asymptotic behavior of critical points of $\eqref{Eeps}$ defined on a closed Riemannian manifold of dimensions $m\ge 2$ has been performed by Lin and Wang in \cite{LinWa2}.\par
Critical points of \eqref{Eeps} are weak solutions of 
\begin{equation}\label{GLeq}
 -\Delta_h u=u\ \frac{1-|u|^2}{\varepsilon^2}.
 \end{equation}
 Given a sequence $u_{\varepsilon_k}$ of weak solutions to \eqref{GLeq} we expect to find in the (weak) limit a harmonic map $u:\Sigma\to {S}^{n}$. As usual for conformally invariant problems, the conformal group leaves the functional $E$ invariant and concentration of energy can a-priori not be excluded. While the Ginzburg-Landau perturbation at fixed scale $E_\varepsilon$ is not conformally invariant, families of critical points $u_\varepsilon$ for $\varepsilon\to 0$ are a-priori also subject to concentration phenomena. In general, the limit of a sequence of critical points $u_{\varepsilon}$ of the energy \eqref{Eeps} will consist of a harmonic map $u_\infty:\Sigma\to {S}^{n}$ as well as a finite number of bubbles $v_{\infty}^i\colon {S}^2\to {S}^{n}$. The precise statement of the well-known bubble-tree convergence is Theorem \ref{bubbletreeconv}, (see Theorem A in \cite{LinWa}):

\begin{Th*}\label{bubbletreeconv}
 Let $u_k$ be a sequence of critical points of $E_{\varepsilon_k}$, $\varepsilon_k\to 0$, with uniformly bounded energy
 \begin{equation}\label{unifbound}
 \sup\limits_k E_{\varepsilon_k}(u_k)=:\Lambda<\infty.
 \end{equation}
 Then there exist sequences $x_k^i\in\Sigma\to a_i\in \Sigma$, $0<\delta_k^i\to 0$ as $k\to +\infty$, $i\in\{1,\ldots N\}$, harmonic maps $u_\infty\colon\Sigma\to {S}^{n-1}$ and $v_\infty^i \colon {\mathbb{C}}\to {S}^{n}$ such that
 \begin{equation}
 u_k\to u_\infty \text{ in } C^\infty(\Sigma\setminus\{a_1,\ldots a_N\}\})
 \end{equation}
 \begin{equation}\label{energyquant}
 \lim_{k\rightarrow +\infty}\left\|e_{\epsilon_k}(u_k)-e(u_\infty)- e(v_\infty^i)\left(\frac{x-x_k^i}{\delta_k^i}\right)\right\|_{L^{1}_{loc}(\Sigma)}=0.
 \end{equation}
\end{Th*}

 The two main contributions of this paper are the upper continuity of the Morse index and the point-wise control of the density in the neck regions:

 \begin{Th*}\label{theoremneckcontrolc0intr}
Let $u_k:(\Sigma,h)\to\mathbb{R}^{n+1}$ be a sequence of critical points of ${E}_{\varepsilon_k}$ with uniformly bounded energy ${E}_{\varepsilon_k}(u_k)\leq \Lambda$ and assume that it converges in the bubble tree sense to $(u_\infty,v_\infty)$. Then there exist $C>0$, $0<\beta<1$ and $\eta_0>0$, such that for all $\eta<\eta_0$ and for all $x\in A(\eta,\delta_k)=B(0,\eta)\setminus \bar{B}(0,\frac{\delta_k}{\eta})$
 \begin{equation}\label{pointwiseestimateeqintro}
 |x|^2|e_{\varepsilon_k}(u_k)|(x)\leq C\bigg(\frac{|x|^\beta}{\eta^\beta}+\frac{\delta_k^\beta}{\eta^\beta|x|^\beta}\bigg)o_{\eta,k}(1).
 \end{equation}
 where $\lim\limits_{\eta\to 0}\limsup\limits_{k\to +\infty}o_{\eta,k}(1)=0$.
\end{Th*}

\begin{Th*}\label{morseindextheoremfordirichletintr}
Let $u_k:(\Sigma,h)\to\mathbb{R}^{n+1}$ be a sequence of critical points of $ E_{\varepsilon_k}$ with uniformly bounded energy $ E_{\varepsilon_k}(u_k)\leq \Lambda$ and assume that it converges in the bubble tree sense to $(u_\infty,v_\infty)$. Then for $k$ large enough,
\begin{equation}
 \text{Ind}_{ E_{\varepsilon_k}}(u_k)+ \text{Null}_{ E_{\varepsilon_k}}(u_k)\leq \text{Ind}_{ E}(u_\infty)+ \text{Null}_{ E}(u_\infty)+\text{Ind}_{ E}(v_\infty)+\text{Null}_{ E}(v_\infty)\,.
\end{equation}
\end{Th*}
\vspace{0.3cm}
\noindent It should be noted that the lower semi-continuity of the index,
\begin{equation}\label{lowersc}
 \text{Ind}_{ E_{\varepsilon_k}}(u_k)\geq \text{Ind}_{{E}}(u_\infty)+\text{Ind}_{E}(v_\infty).
\end{equation}
is well known in the literature and can be shown using a standard $\log$-cutoff method. The lower-semi-continuity of the index \eqref{lowersc} holds for critical points of \eqref{EulerEq} in any dimension and has been recently used in \cite{KaSt2} by M. Karpukhin and D. Stern to construct non-trivial stationary harmonic maps from manifolds of dimension at least three.\\ \noindent

As a result of our refined neck region analysis, we are also able to show the $L^{2,1}$ energy identity:
\begin{Co*}\label{corollaryl21controlintr}
 Let $u_k:(\Sigma,h)\to \mathbb{R}^{n+1}$ be a sequence of critical points of ${E}_{\varepsilon_k}$ with uniformly bounded energy ${E}_{\varepsilon_k}(u_k)\leq \Lambda$ and assume that it converges in the bubble tree sense to $(u_\infty,v_\infty)$. Then 
 \[\limsup\limits_{k\to\infty}||e_{\epsilon_k}^{1/2}(u_k)||_{L^{2,1}(A(\eta,\delta_k))}=0\,.\]
\end{Co*}

\vspace{0.3cm}

\textbf{Discussion.} The strategy of the proof of theorem \ref{morseindextheoremfordirichlet} is strongly inspired by the authors' and T. Rivière's recent work on the Morse index upper continuity for conformally invariant Lagrangians \cite{DLGR}. The main idea consists in showing that the neck regions of the given sequence \underline{positively} contribute to the second derivative of the Lagrangian $E_{\varepsilon_k}(u_k)$: 
 \begin{equation}\label{secondderivativeint}
 D^2{ E_{\varepsilon_k}}(u_k)[v]=\int_{\Sigma}|\nabla v|^2+\frac{1}{\varepsilon_k^2}\langle u_k,v\rangle^2-\frac{1}{2\varepsilon_k^2}(1-|u_k|^2)|v|^2
 \end{equation}
The key step towards obtaining quantitative estimates for the contribution of the index in the neck regions in our previous work \cite{DLGR} was a precise point-wise estimate of the gradients in these regions (Proposition III.1 in \cite{DLGR}). Similarly in the current paper point-wise estimates in the neck regions for the energy densities
 \begin{equation}\label{density1}
e_{\varepsilon}(u)=\frac{1}{2}|d u|_h^2+\frac{1}{4\varepsilon^2}(1-|u|^2)^2
\end{equation}
play a fundamental role. 
The proof of the point-wise estimates in Theorem \ref{theoremneckcontrolc0intr} requires a subtle analysis of the conservation laws hidden in the Euler-Lagrange equations \eqref{EulerEq}. The approach differs substantially from the one taken in \cite{DLGR} in this part of the proof. More precisely, we first prove that the part of the gradient given by $u_k\wedge\nabla u_k$ satisfies $L^{2,1}$ estimates thanks to the conservation law ${\rm div}(u_k\wedge\nabla u_k)=0$ and then use these estimates to get an improved bound for the entire energy density $e_{\varepsilon_k}(u_k)$ in the neck regions. Unlike in the joint paper \cite{DLGR},
we do not need an a-priori $L^{2,1}$ estimate of the gradient to get \ref{theoremneckcontrolc0}. This is due to the fact that in the case of equation \eqref{GLeq} the conservation law ${\rm div}(u_k\wedge\nabla u_k)=0$ holds on the entire space and we do not need a {\em gauge change} on the neck regions under small assumption of the $L^2$ norm of gradients of the given sequence. As explained in the paper \cite{DLGR} such a change of the Gauge produces a conservation law up to a harmonic error that one needs to control by means of $L^{2,1}$ quantization of the gradient. We would like to mention that the $L^2$-energy identity and the no-neck property for a sequence of critical points of the energy \eqref{Eeps} with uniformly bounded energy have been already obtained by different methods by F. H. Lin and C. Wang in \cite{LinWa2}. Our analysis yields a precise {\em quantitative} version of the no-neck property, which is necessary for the proof of Theorem \ref{morseindextheoremfordirichletintr}.\par
 We believe that our approach will be robust enough to treat the case of Ginzburg - Landau type approximation of conformally invariant variational problems in $2$-D of the form
 \be
\label{GGLE}
\mathfrak{L}_\varepsilon(u):=\int_\Sigma \frac{1}{2}|du|^2_h\ dvol_h +\frac{1}{\varepsilon^2}F(u)dvol_h+u^*\alpha
\ee 
 
 where $F$ is some smoothing of the square of the distance to a $C^2$ closed Riemannian manifold ${\mathcal{N}}^n$ of arbitrary dimension $n$, and $\al$ is an arbitrary $2$-form on ${\mathcal{N}}^n$. 

 The analysis of the extended index stability for the Lagrangians \eqref{GGLE} is the aim of our future investigation.\par

\vspace{0.5cm}
\textbf{Further directions. }We finally mention that M. Karpukhin and D. Stern in \cite{KaSt} characterized maximal conformal eigenvalues of surfaces using the energy of harmonic maps coming from a converging sequence of Ginzburg-Landau critical points taking values in $S^n$, with an explicit upper bound for the indices of the critical points along the sequence. The bubbling of (some of) the min-max sequences constructed in \cite{KaSt} could be ruled out by the authors using index considerations. It would be interesting to understand which consequences and/or applications the upper-semi-continuity of the index for this sequences as proven in Theorem \ref{morseindextheoremfordirichlet} might have.

\vspace{0.5cm}
\textbf{Acknowledgements. }The authors are grateful to Prof. Tristan Rivi\`ere for many useful discussions on the subject as well as for his valuable comments and advice.

\setcounter{section}{0}
\setcounter{subsection}{0}

\section{Preliminaries and notation}
We consider an arbitrary smooth closed and oriented 2-dimensional Riemann manifold $(\Sigma,h)$ and 
the Ginzburg-Landau perturbation of the Dirichlet energy given by
\begin{equation}\label{GLenergy}
 E_\varepsilon(u)= \int_{\Sigma}e_\varepsilon(u)\,dvol_h
\end{equation}
 where $u\in W^{1,2}(\Sigma,\mathbb{R}^{n+1})$ and 
 \begin{equation}\label{energydens}
 e_{\varepsilon}(u)=\frac{1}{2}|d u|_h^2+\frac{1}{4\varepsilon^2}(1-|u|^2)^2\,.
 \end{equation}
Critical points of the functional \eqref{GLenergy} satisfy the Euler-Lagrange equation
\begin{equation}\label{EEeps}
-\Delta_h u=u\ \frac{1-|u|^2}{\varepsilon^2}~~\mbox{in $\Sigma$}.
\end{equation}
By standard elliptic regularity arguments such critical points are smooth, and the following uniform $L^\infty$ bounds hold:
\begin{Lm}\label{supestimatelemma}
 Let $u\colon (\Sigma,h)\to\mathbb{R}^{n+1}$ be a smooth solution of 
 \begin{equation}\label{pde}
 -\Delta_h u = u\ \frac{1-|u|^2}{\varepsilon^2}
 \end{equation}
 for some $\varepsilon>0$. Then
 \begin{equation}\label{linftyest}
 \|u\|_{L^\infty(\Sigma)}\leq 1.
 \end{equation}
\end{Lm}
\noindent {\bf Proof of Lemma \ref{supestimatelemma}.}
 We recall the following identity
\begin{equation}
 \frac{1}{2}\Delta_h |u|^2
= {\rm div}_h(u\cdot du)=u\Delta_h u+|d u|_h^2\,.
\end{equation}
Multiplying \eqref{pde} by $u(1-|u|^2)^-$, where $f^-=\min(0,f)\leq 0$ and integrating over $\Sigma$ we get
\begin{align*}
 \int_\Sigma-(1-|u|^2)^-u\Delta_h u=\int_\Sigma -(1-|u|^2)^-\bigg(\frac{1}{2}\Delta_h |u|^2-|d u|_h^2\bigg)=\int_\Sigma\frac{|u|^2}{\varepsilon^2}\bigg((1-|u|^2)^-\bigg)^2.
\end{align*}
This implies
\begin{align}
 &\int_\Sigma - d (1-|u|^2)^-\cdot_h \bigg(\frac{d |u|^2}{2}\bigg)\,dvol_h=\int_\Sigma (1-|u|^2)^-\frac{1}{2}\Delta_h |u|^2 
 \\ \noindent=&\int_\Sigma (1-|u|^2)^-|d u|_h^2-\int_\Sigma\frac{|u|^2}{\varepsilon^2}\bigg((1-|u|^2)^-\bigg)^2 \,dvol_h
\leq 0\,.\nonumber
\end{align}
On the other hand we have 
\begin{equation}
 -\frac{1}{2} d(1-|u|^2)^-\cdot_h (-d(1-|u|^2))=\frac{1}{2}|d(1-|u|^2)^-|^2 \,dvol_h\geq 0.
\end{equation}
Hence $|u|\le 1$ a.e in $\Sigma$ and we conclude the proof.
\qed

\vspace{0.8cm}

\noindent Let $\varepsilon_k$ be a monotone sequence converging to zero and consider an associated family of critical points $u_k$ of $E_{\varepsilon_k}$, 
\begin{equation}
 -\Delta_h u_k=u_k\frac{1-|u_k|^2}{\varepsilon_k^2}
\end{equation}
with uniformly bounded energy in the sense
\begin{equation}
 E_{\varepsilon_k}(u_k)=\int_{\Sigma}\frac{1}{2}|\nabla u_k|_h^2+\frac{(1-|u_k|^2)^2}{4\varepsilon_k^2}\,dvol_h\leq\Lambda<\infty\,.
\end{equation}
Let us denote by $u_\infty$ the weak limit
\begin{equation}
 u_k\rightharpoonup u_\infty\text{ weakly in }W^{1,2}(\Sigma,\mathbb{R}^{n+1}).
\end{equation}

 \noindent While each $u_k$ is smooth, no uniform estimates for higher derivatives a priori hold, so that the weak convergence to $u_\infty$ cannot be expected to be strong in general. Thanks to the so-called small-energy regularity, one can prove that the convergence is strong away from a finite number of points, and that around these points, the convergence is strong after a re-scaling. This kind of convergence is known as bubble-tree convergence and is common for conformally invariant problems.
 By following the approach of R. Schoen in \cite{Sch} for the case of harmonic maps we get the following the small-energy regularity result whose proof is postponed in the appendix.
\begin{Prop}\label{propereg}
Let $u:B_r\to\mathbb{R}^{n+1}$ satisfy
\begin{equation}
 -\Delta u=u\ \frac{1-|u|^2}{\varepsilon^2}
\end{equation}
for some $\varepsilon>0$. There exists $\delta_0$ such that if
\begin{equation}\label{delta0}
 \int_{B_r}\frac{1}{2}|\nabla u|^2+\frac{(1-|u|^2)^2}{4\varepsilon^2}<\delta_0
\end{equation}
then
\begin{equation}\label{Linfty}
 r^2\bigg[\|\nabla u||^2_{L^{\infty}(B_{r/2})}+\frac{1}{\varepsilon^2}\big\vert\big\vert(1-|u|^2)^2\big\vert\big\vert_{L^\infty(B_{r/2})}\bigg]\leq C\int_{B_r}e_\varepsilon(u)
\end{equation}
for some $C>0$ independent of $\varepsilon.$
 \end{Prop}
\noindent From Proposition \ref{propereg} and Lemma \ref{supestimatelemma}, as observed by Lin and Wang in \cite{LinWa2} (Lemma 2.2), one can obtain second order point-wise estimates.

\begin{Co}\label{deltareglap}
There exists $0<\delta_0<1$, such that for any critical point $u$ of $E_\varepsilon$ satisfying
 \begin{equation}
 \int_{B_r}e_\varepsilon(u)<\delta_0,
 \end{equation}
 for some $r>0$, then
 \begin{equation}
r^2\|\Delta u\|_{L^\infty(B_{r/4})}\leq C{\displaystyle{\int_{B_r}e_\varepsilon(u)}} +\frac{Cr^2}{\epsilon^2}e^{-\frac{Cr}{\varepsilon}}\,.
 \end{equation}
\end{Co}
\noindent {\bf Proof of Corollary \ref{deltareglap}.}
Since the estimate is scale-invariant, we can assume $r=2$. We have:
\begin{equation}
 -\Delta(1-|u|^2)=2(u\Delta u+|\nabla u|^2)=-\frac{2|u|^2(1-|u|^2)}{\varepsilon^2} +2|\nabla u|^2
\end{equation}
The following estimate holds
\begin{eqnarray}\label{laplest}
 -\varepsilon^2\Delta(1-|u|^2)+2 (1-|u|^2)&= &-{2|u|^2(1-|u|^2)} +2\varepsilon^2|\nabla u|^2+2 (1-|u|^2)\nonumber\\ \noindent
 &\le& C\varepsilon^2 \int_{B_2}e_\varepsilon(u).
 \end{eqnarray}
 where in the last row of \eqref{laplest} we have used the fact that $(1-|u|^2)\le 1$ and applied the estimate \eqref{Linfty} in $B_1$

We observe that for all $\varepsilon\le 1$, $f=e^{\frac{1}{2\varepsilon}(|x|^2-1)}$ is a super-solution of $-\varepsilon^2\Delta f+2f=0$ for $x\in B_1$.
Then $w=1-|u|^2-f$ satisfies
\begin{equation}
\left\{\begin{array}{cc}
-\varepsilon^2\Delta w+2 w\le 0& ~~\mbox{in $B_1$}\\ \noindent
w\le 0
\end{array}\right.\end{equation}
Maximum principle yields 
\begin{equation}
w(x)\le C\varepsilon^2 \int_{B_2}e_\varepsilon(u)\,.
\end{equation}
Hence
\begin{equation}
 \frac{1-|u|^2}{\varepsilon^2}\leq \frac{1}{\epsilon^2}e^{\frac{1}{2\varepsilon}(|x|^2-1)}+C\int_{B_2}e_\varepsilon(u)\,.
\end{equation}
By re-scaling one obtains the statement at arbitrary scale $r$. This concludes the proof of the Corollary \ref{deltareglap}.
\qed
\vspace{0.8cm}\\ \noindent
 \noindent For the reader's convenience we propose here {\bf an alternative proof of Theorem \ref{bubbletreeconv}} to that obtained by Lin and Wang in \cite{LinWa}. This also gives us the occasion to introduce the notation that will be used in the later chapters.\par\vspace{0.3cm}
\noindent {\bf Proof of Theorem \ref{bubbletreeconv}.} We divide the proof into several steps.\par\vspace{0.5cm}
 \noindent{\bf Step 1: weak limits.} Since $\|u_k\|_{W^{1,2}}\leq CE_{\epsilon_k}(u_k)\leq\Lambda$, $u_k\rightharpoonup u$ weakly in $W^{1,2}(\Sigma,\mathbb{R}^{n+1})$ {\bf Claim:} $u_{\infty}\in W^{1,2}(\Sigma,S^n)$. \par
\noindent {\bf Proof of the claim.} Assume by contradiction there exists $L\subset\Sigma$ of positive measure $\mu(L)>0$ such that $\forall x\in L, |u_\infty(x)|<\sqrt{1-\gamma}$ for a fixed $\gamma\in (0,1)$. By Egorov's Theorem, $u_k$ converges uniformly outside of a measurable set $N$ of arbitrarily small measure so that $u_\infty$ is continuous on $\Sigma\setminus N$. Therefore for $k$ large enough we have $ |u_k|^2<1-\displaystyle\frac{\gamma}{2} \text{ on } L\setminus N$. But then ${{E}}_{\varepsilon_k}(u_k)\geq \mu(L\setminus N)\frac{\gamma^2}{16\varepsilon_k^2}$ which contradicts the bounded energy assumption. \par\vspace{0.3cm}
\noindent{\bf Step 2: concentration points and detection of bubbles.}\footnote{
 A {\bf Bubble} is a {\bf non-constant} harmonic map $u\in \dot W^{1/2}(\mathbb{C},{S}^{m-1})$\,. 
}

\noindent By Proposition \ref{propereg}, the sequence $u_k$ strongly converges to $u_\infty$ in $C^{\ell}_{loc}(\Sigma\setminus\{a_1,\ldots,a_N,\R^m)\}$ for all $\ell\in\N$.
 To simplify the presentation, we suppose that $N=1$.

\noindent Let $\delta_0>0$ be the constant coming from Proposition \ref{propereg}, and define
\begin{equation}\label{weakc4}
 \varrho_k=\inf\bigg\{ \varrho>0\,\,\big\vert\,\, \mbox{there is}~x\in\Sigma\,, ~\mbox{such that}~ B(x, \varrho)\subset\Sigma,~\mbox{and}~\int_{B_{\varrho(x)}} e_{\varepsilon_k}(u_k)=\frac{\delta_0}{2}\bigg\}
\end{equation}
 
\noindent {\bf Case 1.} Up to subsequence $\varrho_k\to\varrho>0$ as $k\to+\infty$. In this case there is no concentration and 
 $$u_k\to u_\infty~~~\mbox{as $k\to+\infty$ in $C^{\ell}_{loc}(\Sigma)$}$$
 {\bf Case 2.} Up to subsequence $\varrho_k\to 0$ as $k\to+\infty$.
 Let $x_k\in\Sigma$ be such that
\begin{equation}
 \int_{B_{\varrho_k}(x_k)}e_{\varepsilon_k}(u_k)=\frac{\delta_0}{2}.
\end{equation}
By compactness of $\Sigma$, $x_k$ converges (up to subsequence) to $a_1$ as $k\to +\infty$ (outside every neighborhood of $a_1$ there is no concentration of the energy).\par
\noindent We rescale the maps $u_k$ around the concentration point $x_k$\footnote{here in conformal coordinates, alternatively one could also use the exponential map}
\begin{equation}
 v_k(z)=u_k(\varrho_k z+x_k)\,,
\end{equation}
 $v_k$ satisfies the equation
\begin{equation}\label{equationbubble}
 - \Delta v_k = \bigg(\frac{\varrho_k}{\varepsilon_k}\bigg)^2(1-|v_k|^2)v_k\,.
\end{equation}
Assume for simplicity that there is only one such sequence $x_k\to a\in A$. Then, for all $y\in\mathbb{C}$ and $k$ large enough
\begin{equation}\label{energyofblowps}
 \frac{\delta_0}{2}\geq \int_{{ {B(\varrho_k y+x_k,\varrho_k)}}}e_{\varepsilon_k}(u_k)\,dvol_h=\int_{B(y,1)}e_{\varepsilon_k/\varrho_k}(v_k) \,dvol_h\,,
\end{equation}
which shows that the blow-ups do not have concentration points and therefore
\begin{equation}\label{smoothconvergencebubble}
 v_k\xrightarrow{C^\infty_{loc}}v_\infty\,.
\end{equation}
In particular $\|\nabla v_\infty\|_{L^2}\leq \liminf\limits_k\|\nabla v_k\|_{L^2}\leq C\Lambda$ i.e. $v_\infty\in {\dot{W}}^{1,2}(\mathbb{C},S^n)$.\par\vspace{0.3cm}
\noindent\underline{Claim:} The limit $v_\infty$ has a holomorphic Hopf differential.\\ \noindent
By smooth convergence \eqref{smoothconvergencebubble},
\begin{equation}\label{equationatbubblelimit}
 - \lim\limits_k\Delta v_k = \lim\limits_k\bigg(\frac{\varrho_k}{\varepsilon_k}\bigg)^2(1-|v_k|^2)v_k=-\Delta v_\infty
\end{equation}
and therefore for any $z\in\mathbb{C}$
\begin{equation}
 \lim\limits_k \bigg(\frac{\varrho_k}{\varepsilon_k}\bigg)^2(1-|v_k(z)|^2)v_k(z)<\infty.
\end{equation}
On the other hand, the pointwise convergence of the $v_k$'s guarantees that for any fixed $z\in \mathbb{C}$, there exists $N_z\in\mathbb{N}$ large enough so that $|v_k-v|\leq \frac{1}{2}$ for all $k\geq N_z$ which combined with $|v_\infty(z)|=1$ yields 
\begin{equation}\label{lowerboundinlimit}
 \inf\limits_{k\geq N_z}|v_k(z)|\geq\frac{1}{2}.
\end{equation}
Combining \eqref{equationatbubblelimit} and \eqref{lowerboundinlimit} gives
\begin{equation}
 2\lim\limits_k \bigg(\frac{\varrho_k}{\varepsilon_k}\bigg)^2(1-|v_k(x)|^2)\leq \lim\limits_k \bigg(\frac{\varrho_k}{\varepsilon_k}\bigg)^2(1-|v_k(z)|^2)|v_k(z)| =|\Delta v_\infty(z)|<\infty
\end{equation}
so that the factor
\begin{equation}\label{convergenceviscosityterm}
 \bigg(\frac{\varrho_k}{\varepsilon_k}\bigg)^2(1-|v_k(z)|^2)\rightarrow g(z)
\end{equation} converges pointwise. The Hopf differential of the $v_k$'s is given by
\begin{equation}\label{Hopfk}
 H_k(z)= \partial_z v_k(z)\cdot \partial_z v_k(z)
\end{equation}
and by differentiating \eqref{Hopfk} with respect to $\overline{z}$ we get
\begin{equation}\label{Hopfk2}
 \partial_{\overline{z}}H_k(z)=2\partial_z v_k\cdot \Delta v_k=-2\partial_z v_k\cdot \bigg(\bigg(\frac{\varrho_k}{\varepsilon_n}\bigg)^2(1-|v_k|^2)v_k\bigg).
\end{equation}
By locally smooth convergence \eqref{smoothconvergencebubble}, for any $z\in\mathbb{C}$,
\begin{equation}
 H_k(z)\rightarrow H_\infty(z)
\end{equation}
and by \eqref{convergenceviscosityterm}
\begin{equation}
 \partial_{\overline{z}}H_\infty(z)= -2\partial_zv_\infty(z)\cdot g_\infty(z)v_\infty(z)\,.
\end{equation}
Since $v_\infty$ is smooth and takes values in the sphere, $\nabla v_\infty\perp v_\infty$ implies $\partial_{\overline{z}}H_\infty=0$. This implies that $v\colon \mathbb{C}\to S^{n}$ is a weakly harmonic map, namely $v_\infty\wedge\Delta v_\infty=0$.
\footnote{one could also use smooth convergence to pass to the limit in the equation $v_k\wedge\Delta v_k=0$ to prove the claim}
We also notice that from $\nabla v_\infty\in L^2$ we deduce $H_\infty\in L^1$, so $H_\infty:\mathbb{C}\rightarrow\mathbb{C}$ is a holomorphic integrable function and therefore
\begin{equation}\label{zerohopfdifferential}
 H_\infty\equiv 0\,.
\end{equation}
 By the regularity theory and point removability property for weakly harmonic maps in two dimensions, $v_\infty$ can be identified with a smooth harmonic map
$$
 v_\infty \colon S^2\to S^{n}.$$
By classical results on harmonic maps there exists a constant\footnote{actually independent of $n$, see \cite{KaSt}} such that
\begin{equation}\label{lowerb}
 E(v_\infty)\geq C(n)>0.
\end{equation}
From the lower bound \eqref{lowerb} we deduce the existence of at most finitely many {\em bubbles} around $a_1$. \par\vspace{0.3cm}
\noindent\textbf{Step 3: $L^2$ energy quantization.} 
 From the discussion in step 2, up to a subsequence, we have the following weak convergence in the sense of Radon measures (we have assumed that there is one concentration point):
\begin{equation}\label{decommeasures}
 e_{\epsilon_k}(u_k)\,dvol_h\rightharpoonup \frac{1}{2}|du_\infty|_h^2\,dvol_h+\nu\delta_{a_1}
\end{equation}
as $k\to+\infty$, with $\nu\ge 0\,.$ \par
 We also assume for simplicity that \underline{there is only one bubble}.
 \par

 Let $\phi\in C(\Sigma)$, we have:
\begin{equation}\label{neck1}
 \int_{B(x_k,\varepsilon_k)}\phi|d u_k|_h^2\,dvol_h=\int_{B(0,{\varepsilon_k/\varrho_k})}\phi(\varrho_kx+x_k)|d u_k(\varrho_k x+x_k)|_h^2\,dvol_h
\end{equation}
and choosing the sequence $\varepsilon_k=\eta^{-1}\varrho_k$, yields
 \begin{align}\label{convergencemicroscopicscale}
 \lim\limits_{k\to\infty}\int_{B(x_k,\varepsilon_k)}\phi |d u_k|_h^2\,dvol_h&=\lim\limits_{k\to\infty}\int_{B(0,{\eta^{-1})}}\phi(\varrho_k x+x_k)
 |du_k(\varrho_k x+x_k)|_h^2\,dvol_h\nonumber\\ \noindent
 &=\lim\limits_{k\to\infty}\int_{B(0,{\eta^{-1})}}\phi(\varrho_k x+x_k)
 |dv_k(x)|_h^2\,dvol_h\nonumber\\ \noindent
 &=\phi(a) \int_{B(0,{\eta^{-1})}}|dv_{\infty}(y)|_h^2\,dvol_h\,. \end{align}
 On the other hand for all $\eta>0$ we have 
 \begin{equation}\label{convergencemacroscopicscale}
 \lim\limits_{k\to\infty}\int_{\Sigma\setminus B(x_k,\eta)}\phi|d u_k|_h^2\,dvol_h= \int_{\Sigma\setminus B(a,\eta)}\phi|d u_\infty|_h^2\,dvol_h\end{equation}
 The computation \eqref{convergencemicroscopicscale} and \eqref{convergencemacroscopicscale} show that
 \begin{equation}\label{ineqenergyquant}
 |d u_\infty|_h^2\,dvol_h+\nu \delta_{a}\geq |d u_\infty|_h^2\,dvol_h+2E(v_\infty) \delta_{a}\,.
 \end{equation}
 To conclude the proof of the Theorem it remains to show that \eqref{ineqenergyquant} is actually an equality.\\ \noindent

\noindent Combining \eqref{convergencemicroscopicscale}, \eqref{convergencemacroscopicscale} we can decompose
\begin{equation}\label{ggg}
\int_{\Sigma}\phi e_{\epsilon_k(u_k)} = \underbrace{\int_{B(x_k,\eta^{-1}\delta_k)}\phi e_{\epsilon_k(u_k)}\,dvol_h}_{\to \phi(a_1) \int_{B(0,{\eta^{-1})}}\frac{1}{2}|dv_{\infty}(y)|_h^2\,dvol_h}+ \underbrace{\int_{\Sigma\setminus B(x_k,\eta)}\phi e_{\epsilon_k(u_k)}}_{\to \int_{\Sigma}(1-\11_{B(x_k,\eta)})\phi \frac{1}{2}|d u_\infty|_h^2} + {\int_{B(x_k,\eta)\setminus B(x_k\eta^{-1}\delta_k)}\phi e_{\epsilon_k}(u_k)}.
\end{equation}
We also observe that by dominated convergence we have
\begin{equation}\label{gggg}
 \lim\limits_{\eta\to 0} \int_{\Sigma}(1-\11_{B(x_k,\eta)})\phi|d u_\infty|_h^2\,dvol_h= \int_{\Sigma}\phi|d u_\infty|_h^2\,dvol_h.
\end{equation}
In \cite{LinWa2} (Lemma 3.1) Lin and Wang\footnote{The proof relies on the super-convexity of the tangential energy. The method of Lin and Wang provides decay rates for $\int_{B(x_k,\eta)\setminus B(x_k,\eta^{-1}\delta_k)}|\partial_\theta u_k|^2\,dvol_h$, it is however a-priori not clear whether these decay rates still hold for the entire energy density $e_{\epsilon_k}(u_k)$. In order to prove our main theorem, it is essential for us to know this behaviour. The tools introduced in section $3$ will eventually allow us to obtain estimate \eqref{estekbeta}, which is a quantitative version of the no-neck energy.} showed the following no-neck energy result:
\begin{equation}\label{noneckclaim}
 \lim\limits_{\eta\to 0}\limsup\limits_{k\to\infty} \int_{B(x_k,\eta)\setminus B(x_k,\eta^{-1}\varrho_k)}e_{\epsilon_k}(u_k)\,dvol_h=0\,.
\end{equation}
Combining \eqref{ggg}, \eqref{gggg} and \eqref{noneckclaim} we get,
\begin{equation}
\lim\limits_{\eta\to 0}\limsup\limits_{k\to\infty}\int_{\Sigma}\phi|d u_k|_h^2\,dvol_h=\int_{\Sigma}\phi|d u_{\infty}|_h^2\,dvol_h+\phi(a)\int_{S^2}|dv_{\infty}(y)|_h^2\,dvol_h\,,
\end{equation}
which implies 
\begin{equation}\label{eqenergyquant}
 |d u_\infty|_h^2\,dvol_h+\nu \delta_{a_1} = |d u_\infty|_h^2\,dvol_h+2E(v_\infty) \delta_{a_1}\,
 \end{equation}
 The proof of \eqref{noneckclaim} is given. Hence we can conclude the proof of step 2 and of Theorem \ref{bubbletreeconv}.
 \qed
 \section{Point-wise estimates of the energy density in the neck regions}
 
The goal of this section is to use the tools developed by the authors in collaboration with Rivi\`ere in \cite{DLGR} to obtain point-wise estimates of the energy density in the neck regions. It is important to have a point-wise control of the entire energy density and not just for the gradient, since this would not be enough to effectively control the second variation. The methods in \cite{DLGR} rely on the fact that the critical points of $\mathcal{L}$ enjoy compensation compactness properties coming from Wente-type estimates. In the setting of the Ginzburg-Landau perturbation studied here we show that the quantity $u\wedge \nabla u$ satisfies good compensation compactness properties and can be used to bound the entire energy density in an appropriate sense.
\par

 We start showing some preliminary results.
 \begin{Lm}\label{lemmaepsilondeltaratio}
 Let $\varrho>0$ be defined as in \eqref{weakc4}. Then
 \begin{align}
& \liminf\limits_{k\to\infty} \frac{\varrho_k}{\varepsilon_k}=\infty.
 \end{align}
\end{Lm}
\noindent {\bf Proof of Lemma \ref{lemmaepsilondeltaratio}.} 
 Define as in Theorem \ref{bubbletreeconv} for $x:=\varrho_k y$
 \begin{equation}
 v_k(y)=u_k(\varrho_k y+x_k)
 \end{equation}
 then $v_k$ satisfies
 \begin{equation}
 \Delta v_k+\frac{\varrho_k^2}{\varepsilon_k^2}v_k(1-|v_k|^2)=0\,.
 \end{equation}
 We assume by contradiction that there is a subsequence such that $\frac{\varrho_k}{\varepsilon_k}\to c<\infty$. \\ \noindent
 Denoting by $v$ the strong limit $v_k\to v$ in $C^\infty_{loc}(\mathbb{R}^2)$ and integrating by parts with the test function $x_i\displaystyle\frac{\partial v_k}{\partial x_i}$ (Pohozahev identity) yields
 \begin{equation}\label{eq1}
 \int_{B_R}x_i\frac{\partial v}{\partial x_i}(\Delta v+cv(1-|v|^2))=0
 \end{equation}
 \begin{align}\label{eq2}
 \int_{B_R}x_i\frac{\partial v}{\partial x_i}\Delta v&=\int_{\partial B_R}x_i\frac{\partial v}{\partial x_i}\frac{\partial v}{\partial x_k}\frac{x_k}{R}-\int_{B_R}\delta_{ik}\frac{\partial v}{\partial x_i}\frac{\partial v}{\partial x_k}-\int_{B_R}x_i\frac{\partial}{\partial x_i}\frac{|\nabla v|^2}{2}
 \\ \noindent&=\int_{\partial B_R}R\bigg\vert\frac{\partial v}{\partial r}\bigg\vert^2-\int_{B_R}|\nabla v|^2-\int_{B_R}x_i\frac{\partial }{\partial x_i}\frac{|\nabla v|^2}{2}\\ \noindent&=\int_{\partial B_R}R\bigg\vert\frac{\partial v}{\partial r}\bigg\vert^2-\int_{\partial B_R}R\frac{|\nabla v|^2}{2}\nonumber
 \end{align}
 \begin{align}\label{eq4}
 \int_{B_R}x_i\frac{\partial v}{\partial x_i}v(1-|v|^2)&=\frac{1}{2}\int_{B_R}x_i\frac{\partial |v|^2}{\partial x_i}(1-|v|^2)=-\frac{1}{4}\int_{B_R}x_i\frac{\partial (1-|v|^2)^2}{\partial x_i}\\ \noindent
 &=-\int_{\partial B_R}R\frac{(1-|v|^2)^2}{4}+\frac{1}{2}\int_{B_R}(1-|v|^2)^2.
 \end{align}
 Putting the equations \eqref{eq2} and \eqref{eq4} together yields 
 \begin{equation}\label{eq5}
 c\int_{B_R}(1-|v|^2)^2\leq CR\int_{\partial B_R}[|\nabla v|^2+(1-|v|^2)^2] dx
 \end{equation}
 By Fubini's theorem for all $ j\in\mathbb{N}$ there is $ R_j\in [2^j,2^{j+1}]$ such that
 \begin{equation}
 R_j\int_{\partial B_{R_j}}[|\nabla v|^2+(1-|v|^2)^2]dx\leq C\int_{B_{2^{j+1}}\setminus B_{2^j}}|\nabla v|^2+(1-|v|^2)^2\xrightarrow{j\to\infty}0.
\end{equation}
 We conclude
 \begin{equation}
 \int_{\mathbb{R}^2}(1-|v|^2)^2=0\implies |v|=1\implies \Delta v=0,
 \end{equation}
 but the only harmonic function with bounded energy on the plane is $v\equiv 0$, which contradicts the fact that $v$ has strictly positive energy. This concludes the proof of the Lemma.
 \qed
 \par

\begin{Co}\label{gradbound} There is $\bar k $ such that for all $k>\bar k$ and for all $x\in \Sigma$ we have
 \be\label{estgrad1}
 |\nabla u_k|^2(x)\leq Ce_\epsilon(u_k)(x)\leq \frac{C}{\varrho_k^2}\,.
 \ee
\end{Co}
\noindent {\bf Proof of Corollary \ref{gradbound}.} By definition of $\varrho_k$ in \eqref{weakc4}, for all $x\in\Sigma$
\[\int_{B_{\varrho_k}(x)}e_{\epsilon_k}(u_k)\leq\frac{\delta_0}{2}\,,\]
so that by Proposition \ref{propereg},
\begin{equation}
 |\nabla u_k|^2(x)+\frac{(1-|u_k|^2)^2}{\epsilon_k^2}(x)\leq\frac{C}{\varrho_k^2}\int_{B_{\varrho_k}(x)}e_{\epsilon_k}(u_k)\leq C\frac{\delta_0}{\varrho_k^2}\,. \qed
\end{equation}
 
\begin{Co}\label{corollaryunifconv}
 The convergence to the weak limit is uniform, i.e.
 \begin{equation}
 \|u_k-u_\infty||_{L^\infty(\Sigma)}\to 0.
 \end{equation}
\end{Co}
\noindent {\bf Proof of Corollary \ref{corollaryunifconv}.}
Fix $\varepsilon>0$. We prove the existence of $\bar k>0$, such that for all $k\geq \bar k$,
\begin{equation}\label{linftyconv1}
 \|u_k-u_\infty||_{L^\infty(\Sigma)}\leq \varepsilon.
\end{equation}
 Recall that by Step 1 in the proof of Theorem \ref{bubbletreeconv}, $|u_\infty(x)|=1$ for almost every $x\in\Sigma$. We can now argue by contradiction using Lemma \ref{lemmaepsilondeltaratio}.
 Assume there exists a sequence $x^k\in\Sigma$ such that 
 \begin{equation}
 |u(x^k)|< 1-\varepsilon
 \end{equation}
 then by Corollary \ref{gradbound}
 \begin{equation}\label{linftyconv2}
 \forall x\in B_{c\varrho_k}(x^k),\quad |u|\leq 1-\frac{3}{4}\varepsilon
 \end{equation}
 for a constant $c=c(\varepsilon)$ independent of $k$. In particular,
 \begin{equation}\label{linftyconv3}
 \int_{B_{c\varrho_k}}\frac{(1-|u_k|^2)^2}{\varepsilon_k^2}\geq c(\varepsilon)\frac{\varrho_k^2}{\varepsilon_k^2}.
 \end{equation}
 But we are assuming that the sequence $u_k$ has uniform bounded energy 
 ${E}_{\varepsilon_k}(u_k)\leq\Lambda<\infty$, therefore
 \begin{equation}\label{linftyconv4}
 \Lambda\geq C\frac{\varrho_k^2}{\varepsilon_k^2} \to \infty
 \end{equation}
and we get a contradiction. This concludes the proof of Corollary \ref{corollaryunifconv}. \qed \\ \noindent

We observe that thanks to Theorem \ref{bubbletreeconv} and Corollary \ref{corollaryunifconv} we may suppose that for $k$ large enough we have $|u_k|\ge 3/4$ in $\Sigma.$
 
 \subsection{Hodge decomposition of $u_k\wedge \nabla u_k$}\hfill\vspace{0.3cm}\\ \noindent
 Thanks to conservation laws hidden in the Ginzburg-Landau equation we will show that there is a suitable decomposition of $u_k\wedge\nabla u_k$ in terms of Jacobians and harmonic functions so that one can obtain more precise estimates analogous to the ones obtained by the authors in Section 3 of \cite{DLGR}. \par
We introduce some notation.
 For all $j\in\Z$, $0<\delta<\eta$ we set:

\begin{equation}\label{ring}
 A_j=B_{2^{-j}}\setminus B_{2^{-j-1}} \text{ the dyadic annuli }
 \end{equation}
\begin{equation}\label{neck}
 A(\eta,\delta)=B_\eta\setminus B_{\eta/\delta} 
\end{equation}
 
 All the balls in these notations are centered around the concentration point, which for simplicity of notation is assumed to be unique.\\ \noindent

\begin{Lm}\label{lemmahodgedec}
 Let $u_k:(\Sigma,h)\to\mathbb{R}^{n+1}$ be a sequence of critical points of ${ E}_{\varepsilon_k}$ with uniformly bounded energy ${ E}_{\varepsilon_k}(u_k)\leq \Lambda$. and converging in the bubble tree sense to $(u_\infty,v_\infty)$. Then for all $0<\eta<1$ and $0<\delta_k\to 0$ as $k\to Infty$ we have
 \begin{enumerate}
 \item[i)] there exists $\tilde{u}_k:\mathbb{C}\to \mathbb{R}^{n+1}$ such that
 \begin{equation}
 \tilde{u}_k=u_k \text{ on }A(\eta,\delta_k)\text{ and } \int_{\mathbb{C}}|\nabla \tilde{u}_k|^2\leq C\int_{A(\eta,\delta_k)}|\nabla u_k|^2
 \end{equation}
 \item[ii)] there are functions $h_k$ on $A(\eta,\delta_k)$ and $\varphi_k$ such that
 \begin{eqnarray}
 \Delta\varphi_k&=&\nabla^\perp \tilde{u}_k\wedge \nabla\tilde{u}_k~~~\mbox{in $\mathbb{C}$}
\label{wenteequ}\\ \noindent
 \Delta h_k&=& 0~~~\mbox{in $A(\eta,\delta_k)$}
 \end{eqnarray}
 \begin{equation}\label{hod}
 u_k\wedge \nabla u_k = -\nabla^\perp\varphi_k+\nabla h_k\text{ on }A(\eta,\delta_k).
 \end{equation}
 \item[iii)] $||u_k\wedge\nabla u_k||_{L^{2,1}(A(\eta,\delta_k))}\leq C||\nabla u_k||^2_{L^2(A(2\eta,\delta_k))}$\label{l21controlconservationpart}
 \end{enumerate}

\end{Lm}
\noindent {\bf Proof of Lemma \ref{lemmahodgedec}.}
For the proof of the Whitney-type extension $i)$ we refer to the Appendix C of \cite{DLGR}.\\ \noindent
 It follows from equation \eqref{GLeq} that
 \begin{equation}\label{conslaw}
 \text{div}(u_k\wedge \nabla u_k)=0\quad\text{and}\quad\text{curl}(u_k\wedge\nabla u_k)=\nabla^\perp u_k\wedge \nabla u_k\quad\text{ in }\Sigma\,.
 \end{equation}
 Let $\varphi_k$ be the solution of
 \begin{equation} \Delta\varphi_k=\nabla^\perp\tilde{u_k}\wedge\nabla\tilde{u_k}~~~\text{ on }\mathbb{C}\,.
 \end{equation}
 This yields the conservation laws
 \begin{equation}\label{equationconservationlaws}
 \begin{cases}
 \text{div}(u_k\wedge\nabla u_k -\nabla^\perp \varphi_k)=0~~~\text{ on } A(\eta,\delta_k)\\ \noindent
 \text{curl}(u_k\wedge\nabla u_k -\nabla^\perp \varphi_k)=0~~~\text{ on } A(\eta,\delta_k)\,.
 \end{cases} 
 \end{equation}
 By Wente's estimates
 \begin{equation}\label{estnablavarphi}
 \|\nabla \varphi_k\|_{L^{2,1}(\mathbb{C})}\le C\|\nabla \tilde{u}_k\|^2_{L^{2}(\mathbb{C})}.
 \end{equation}
 Now we observe that for every $\rho\in [\eta^{-1}\delta_k, \eta]$ it holds:
 \begin{equation}\label{flow1}
 \int_{\partial B_\rho}(u_k\wedge\nabla u_k-\nabla^{\perp}\varphi_k)\cdot\tau=\int_{\partial B_\rho} \frac{1}{\rho}u_k\wedge\partial_\theta u_k -\partial_r\varphi_k
 \end{equation}
 and
 \begin{eqnarray}
 \int_{\partial B_\rho}\partial_{r}\varphi_k dl&=&\int_{B_\rho}\Delta\varphi_k=\int_{B_\rho}\nabla^\perp\tilde{u}_k\wedge\nabla{u}_k\label{flow2}\\ \noindent
 &=&\int_{B_\rho}\text{div}(\nabla^\perp \tilde{u}_k\wedge \tilde{u}_k)=\int_{\partial B_\rho}\tilde{u}_k\wedge\frac{1}{\rho}\partial_\theta \tilde{u}_k\nonumber
 \end{eqnarray}
By combining \eqref{flow1} and \eqref{flow2} we get that
\begin{equation}
 \int_{\partial B_\rho}(u_k\wedge\nabla u_k-\nabla^{\perp}\varphi_k)\cdot\tau=0
 \end{equation}
 Therefore Hodge's decomposition gives the existence of a potential $h_k$ such that
 \begin{equation}\label{hod2}
 u_k\wedge \nabla u_k - \nabla^\perp \varphi_k=\nabla h_k~~~\text{ on }A(\eta,\delta_k).
 \end{equation}
 which by \eqref{equationconservationlaws} is harmonic. This concludes the proof of the Lemma.\\ \noindent \qed\\ \noindent
 
 \begin{Rm}
 The harmonic functions $h_k$ on the annulus $A(\eta,\delta_k)$ are well understood in terms of their frequencies:

\[
{h}_k={h}_k^++{h}_k^-+{h}_k^0\quad {h}_k^+= \mathrm{Re}\sum_{n>0} h_{n,k}z^n\quad h_k^-= \mathrm{Re}\sum_{n<0} h_{n,k}z^n\quad{h}_k^0={h}_{0,k}+C_{\eta,k}^0\,\log|z|\ .
\]

 \begin{equation}
 h_k=h^0_k+h^+_k+h^-_k=h^+_k+h^-_k+C^0_{\eta,k}\log|x|
 \end{equation}
 We observe that the logarithmic part in this case is zero:
 \begin{equation*}
 \int_{\partial B_\rho}\partial_\nu h_k=2\pi C^0_{\eta,k}=\int_{\partial B_\rho}u_k\wedge\partial_\nu u_k-\frac{1}{\rho}\partial_\theta\varphi_k=\int_{B_\rho}\text{div}(u_k\wedge\nabla u_k)=0\implies C^0_{\eta,k}=0\,.
 \end{equation*}
 \end{Rm}

\subsection{Improved $L^2$-estimates in the neck region}
\begin{Lm}\label{lemmacontrolpotentialpart} Let $u_k$ be a critical point of ${E}_{\varepsilon_k}$ satisfying $\int_{B_R\setminus B_{r}}e_{\varepsilon_k}(u_k)<\delta_0$, for $R/r\in [2,8]$. Then, for every $\sigma>0$,
\begin{align}
 \int_{B_R\setminus B_{r}}|u_k|\frac{(1-|u_k|^2)^2}{\varepsilon_k^2}&\leq \sigma \bigg( E_{\varepsilon_k}(u_k;{B_{2R}\setminus B_{16^{-1}R}}) +\frac{CR^2}{\epsilon_k^2}e^{-\frac{CR}{\varepsilon}}\bigg)^2\\ \noindent&+C\frac{\varepsilon_k^2R^{-2}}{\sigma}\int_{B_R\setminus B_r}\frac{(1-|u_k|^2)^2}{\varepsilon_k^2}\,.
 \end{align}
\end{Lm}
\noindent {\bf Proof of Lemma \ref{lemmacontrolpotentialpart}.}
 \begin{align}\label{L21}
 &\int_{B_R\setminus B_{r}}|u_k|\frac{(1-|u_k|^2)^2}{\varepsilon_k^2}\leq C \int_{B_R\setminus B_{r}}|\Delta u_k|(1-|u_k|)\\ \noindent&\leq C\bigg( E_{\varepsilon_k}(u_k;{B_{2R}\setminus B_{16^{-1}R}}) +\frac{CR^2}{\epsilon_k^2}e^{-\frac{CR}{\varepsilon_k}}\bigg) R^{-2}\int_{B_R\setminus B_r}(1-|u_k|)\\ \noindent
 &\leq C \bigg( E_{\varepsilon_k}(u_k;{B_{2R}\setminus B_{16^{-1}R}}) +\frac{CR^2}{\epsilon_k^2}e^{-\frac{CR}{\varepsilon_k}}\bigg) \bigg(\int_{B_R\setminus B_r}(1-|u_k|)\bigg)^{1/2}R^{-1}\\ \noindent
 &\leq \sigma \bigg( E_{\varepsilon_k}(u_k;{B_{2R}\setminus B_{16^{-1}R}}) +\frac{CR^{2}}{\epsilon_k^2}e^{-\frac{CR}{\varepsilon_k}}\bigg)^2+C\frac{\varepsilon_k^2R^{-2}}{\sigma}\int_{B_R\setminus B_r}\frac{(1-|u_k|^2)^2}{\varepsilon_k^2}\,.
 \end{align}
\qed\\ \noindent

\begin{Lm}\label{lemma15} Let $u_k:\Sigma\to\mathbb{R}^{n+1}$ be a sequence of critical points of ${ E}_{\varepsilon_k}$ satisfying ${ E}_{\varepsilon_k}(u_k)\le \Lambda$ and converging in the bubble tree sense to $(u_\infty,v_\infty)$. Then for all $k$ large enough and for all $j\ge \log_2(\eta/\delta_k)$ there is $\rho_j\in[2^{-j-1},2^{-j}]$ such that for every $\sigma>0$,
 \begin{equation*}
 \int_{B_{2\rho_j}\setminus B_{\rho_j/2}}|\nabla u_k|^2\leq C\int_{B_{2\rho_j}\setminus B_{\rho_j/2}}|u_k\wedge\nabla u_k|^2+C\sigma\int_{\tilde{A}_j}|\nabla u_k|^2+\frac{C\varepsilon_k^2}{\sigma 2^{-2j}}\int_{\tilde{A}_j}\frac{(1-|u_k|^2)^2}{\varepsilon_k^2}
 \end{equation*}
 where $\tilde{A}_j=A_{j-1}\cup A_j\cup A_{j+1}$.
\end{Lm}
\noindent {\bf Proof of Lemma \ref{lemma15}.} Recall that $|u_k|$ satisfies (see \cite{LiRi1}),
 \begin{equation}\label{pdemod}
 \forall x \text{ s.t. } |u_k|(x)\neq0, \quad\Delta|u_k|+|u_k|\bigg(\frac{1-|u_k|^2}{\varepsilon_k^2}\bigg)=\frac{|u_k\wedge\nabla u_k|^2}{|u_k|^3}\,.
\end{equation}
By Corollary \ref{corollaryunifconv} $|u_k|\to 1$ in ${L^{\infty}}$ as $k\to +\infty$, there therefore exists $\bar k$ such that for all $k\ge \bar k$, $|u_k|\ge 3/4$. We first observe that
\begin{equation}\label{decogradient}
|\nabla u_k|^2=|\nabla |u_k||^2+\frac{|\nabla u_k\wedge u_k|^2}{|u_k|^2}
\end{equation}
 
\noindent Multiplying the equation \eqref{pdemod} by $(1-|u_k|)$ and integrating by parts on any $\Omega\subseteq\Sigma$ we get the following estimate
 \begin{eqnarray}\label{estring}
\int_{\Omega}|\nabla u_k|^2&=&\int_{\Omega}\frac{|u_k\wedge\nabla u_k|^2}{|u_k|^2}+|\nabla|u_k||^2\nonumber\\ \noindent&=&\int_{\Omega}\frac{|u_k\wedge\nabla u_k|^2}{|u_k|^2}-\int_{\Omega}\nabla |u_k|\cdot\nabla(1-|u_k|)\nonumber\\ \noindent
&=&\int_{\Omega}\frac{|u_k\wedge\nabla u_k|^2}{|u_k|^2}+\int_{\Omega}\Delta |u_k|(1-|u_k|)-\int_{\partial A_j}(1-|u_k|)\partial_\nu |u_k|\nonumber\\ \noindent
&=&\int_{\Omega}\frac{|u_k\wedge\nabla u_k|^2}{|u_k|^2}\bigg(1+\frac{1-|u_k|}{|u_k|}\bigg)-\int_{\partial \Omega}(1-|u_k|)\partial_\nu |u_k|-\int_{\Omega}\frac{1-|u_k|^2}{\varepsilon_k^2|u_k|^2}(1-|u_k|)\nonumber\\ \noindent
&\leq &\int_{\Omega}\frac{|u_k\wedge\nabla u_k|^2}{|u_k|^3}-\int_{\partial A_j}(1-|u_k|)\partial_\nu |u_k|\,.
\end{eqnarray}
\noindent In order to control the boundary term $\int_{\partial \Omega}(1-|u_k|)\partial_\nu |u_k|$, we observe that the following holds:
\begin{align*}
 \int_{2^{-j-1}}^{2^{-j+1}}\rho\int_{\partial(B_{2\rho}\setminus B_{\rho/2})}|\nabla |u_k||(1-|u_k|)\leq C\bigg(\int_{\tilde{A}_j}|\nabla |u_k||^2\bigg)^{1/2}\bigg(\int_{\tilde{A}_j}(1-|u_k|^2)^2\bigg)^{1/2}
\end{align*}
 For all $j\in\Z$ there must exist $\rho_j\in[2^{-j-1},2^{-j+1}]$ such that
\begin{equation}\label{boundestring2i}
 \rho_j\int_{\partial(B_{2\rho_j}\setminus B_{\rho_j/2})}|\nabla |u_k||(1-|u_k|)\leq C\bigg(\int_{\tilde{A}_j}|\nabla |u_k||^2\bigg)^{1/2}\bigg( \varepsilon_k^2\int_{\tilde{A}_j}\frac{(1-|u_k|^2)^2}{\varepsilon_k^2}\bigg)^{1/2}
\end{equation}
and since $\rho_j\ge 2^{-j-1}$,
\begin{equation}\label{boundestring2ii}
 \int_{\partial(B_{2\rho_j}\setminus B_{\rho_j/2})}|\nabla |u_k||(1-|u_k|)\leq C\bigg(\int_{\tilde{A}_j}|\nabla |u_k||^2\bigg)^{1/2}\bigg(\frac{\varepsilon_k^2}{2^{-2j-2}}\int_{\tilde{A}_j}\frac{(1-|u_k|^2)^2}{\varepsilon_k^2}\bigg)^{1/2}\,,
\end{equation}
\noindent which in turn implies
\begin{equation}\label{boundestring2}
 \int_{\partial {(B_{2\rho_j}\setminus B_{\rho_j/2})}}|\nabla |u_k||(1-|u_k|)\leq \sigma\int_{\tilde{A}_j}|\nabla |u_k||^2 + \frac{C}{\sigma}\frac{\varepsilon_k^2}{2^{-2j-2}}\int_{\tilde{A}_j}\frac{(1-|u_k|^2)^2}{\varepsilon_k^2}\,.
\end{equation}
\noindent By plugging \eqref{boundestring2} into \eqref{estring} in the case $\Omega= {B_{2\rho_j}\setminus B_{\rho_j/2}}$ and using the fact that $|u_k|\ge 3/4$ we get:
\begin{equation*}
 \int_{B_{2\rho_j}\setminus B_{\rho_j/2}}|\nabla u_k|^2\le \int_{B_{2\rho_j}\setminus B_{\rho_j/2}} {|u_k\wedge \nabla u_k|^2}+\sigma\int_{\tilde{A}_j}|\nabla u_k|^2 + \frac{C}{\sigma}\frac{\varepsilon_k^2}{2^{-2j-2}}\int_{\tilde{A}_j}\frac{(1-|u_k|^2)^2}{\varepsilon_k^2}.
\end{equation*}
\qed\\ \noindent
\begin{Co}\label{co1}
By setting $\bar{A}_j=B_{2\rho_j}\setminus B_{\rho_j/2}$, where $\rho_j>0$ comes from Lemma \ref{lemma15}, for $k$ large enough the following holds:
 \begin{equation}
 \int_{\bar{A}_j}|\nabla u_k|^2\leq C\int_{\bar{A}_j}|u_k \wedge \nabla u_k|^2+\sum\limits_{|j'-j|\leq 1}C\sigma \int_{\bar{A}_{j'}}|\nabla u_k|^2+\frac{C\varepsilon_k^2}{\sigma 2^{-2{j'}}}\int_{\bar{A}_{j'}}\frac{(1-|u_k|^2)^2}{\varepsilon_k^2}
 \end{equation}
\end{Co}
\noindent Combining Lemma \ref{lemmacontrolpotentialpart} and Corollary \ref{co1} yields
\begin{Prop}\label{propinneck}
\begin{align}\label{propinneck2}
 \int_{\bar{A}_j}e_{\epsilon_k}(u_k)&\leq C\int_{\bar{A}_j}|u_k \wedge \nabla u_k|^2+\sum\limits_{|j'-j|\leq 3}C\sigma \int_{\bar{A}_{j'}}e_{\epsilon_k}(u_k)+\frac{C\varepsilon_k^2}{\sigma 2^{-2{j'}}}\int_{\bar{A}_{j'}}e_{\epsilon_k}(u_k)\nonumber\\ \noindent+&\sigma \bigg(\frac{C2^{-2j}}{\epsilon_k^2}e^{-\frac{C2^{-j}}{\varepsilon_k}}\bigg)^2
 \end{align}
\end{Prop}
\noindent {\bf Proof of Proposition \ref{propinneck}.} 
By Lemma \ref{lemmacontrolpotentialpart}, which we can apply since the conformal class of $\bar{A}_j$, $2\rho_j/(\rho_j/2)=4$ is bounded,
\begin{align} 
& \int_{\bar{A}_j}\frac{(1-|u_k|^2)^2}{\varepsilon_k^2}\leq \sigma \bigg( E_{\varepsilon_k}(u_k;{B_{4\rho_j}\setminus B_{16^{-1}\rho_j}}) +\frac{C\rho_j^2}{\epsilon_k^2}e^{-\frac{C\rho_j}{\varepsilon_k}}\bigg)^2+C\frac{\varepsilon_k^2\rho_j^{-2}}{\sigma}\int_{\bar{A}_j}\frac{(1-|u_k|^2)^2}{\varepsilon_k^2}\\
&\leq C\sigma \bigg( E_{\varepsilon_k}(u_k;{B_{4\rho_j}\setminus B_{16^{-1}\rho_j}})\bigg)^2 +C\sigma\bigg(\frac{C\rho_j^2}{\epsilon_k^2}e^{-\frac{C\rho_j}{\varepsilon_k}}\bigg)^2+C\frac{\varepsilon_k^2\rho_j^{-2}}{\sigma}\int_{\bar{A}_j}\frac{(1-|u_k|^2)^2}{\varepsilon_k^2}\\
&\leq C\sigma E_{\varepsilon_k}(u_k;{B_{4\rho_j}\setminus B_{16^{-1}\rho_j}}) +C\sigma\bigg(\frac{C\rho_j^2}{\epsilon_k^2}e^{-\frac{C\rho_j}{\varepsilon_k}}\bigg)^2+C\frac{\varepsilon_k^2\rho_j^{-2}}{\sigma}\int_{\bar{A}_j}\frac{(1-|u_k|^2)^2}{\varepsilon_k^2}\\
&\leq C\sigma \sum\limits_{|j'-j|\leq 3}\int_{\bar{A}_{j'}}e_{\epsilon_k}(u_k) +C\sigma\bigg(\frac{C\rho_j^2}{\epsilon_k^2}e^{-\frac{C\rho_j}{\varepsilon_k}}\bigg)^2+C\frac{\varepsilon_k^2\rho_j^{-2}}{\sigma}\int_{\bar{A}_j}\frac{(1-|u_k|^2)^2}{\varepsilon_k^2}\,, \label{moredet}
 \end{align}
 where we used the no-neck property \eqref{noneckclaim} in the second to last step. The Proposition then follows by combining the estimate \eqref{moredet} with Corollary \ref{co1}.\\
 \qed\\ \noindent
 \subsection{Refined analysis in the neck region}\hfill\vspace{0.3cm}\\ \noindent
 \noindent We now combine the Hodge decomposition obtained for $u_k\wedge \nabla u_k$ in \eqref{hod2} with Lemmata F.1, F.2 and G.1 in \cite{DLGR} to obtain an intermediate result for the control of the wedge product.
We will set as above $\bar A_j:=B_{2\rho_j}\setminus B_{\rho_j/2}$. From \eqref{hod2} we get
 \begin{equation}\label{eee1}
 \int_{\bar A_j}|u_k\wedge \nabla u_k|^2\leq 2\int_{\bar A_j}|\nabla h_k|^2+2\int_{\bar A_j}|\nabla\varphi_k|^2. \end{equation}

\begin{Prop}\label{propositionharmpart}
In the setting of Lemma \ref{lemmahodgedec}, there exist $4^{-1}<\mu<1$ and $C>0$ independent of $k$ and $\eta$, such that for $2\frac{\delta_k}{\eta}\leq2^{-s_2}\leq2^{-s_1}\leq\frac{\eta}{2}$,
 \begin{align}\label{esthplusminus}
 &\sum_{j=s_1}^{s_2}\mu^{|j-l|}\int_{\bar A_j}|\nabla h_k^+|^2+ \int_{\bar A_j}|\nabla h_k^-|^2\\ \noindent \leq&\sum_{j=s_1}^{s_2}\mu^{|j-l|}\int_{\tilde{A}_j}|\nabla h_k^+|^2+ \int_{\tilde{A}_j}|\nabla h_k^-|^2
 \\ \noindent\leq &C\bigg(\left(\frac{2^{-l}}{\eta}\right)^{\beta}+\left(\frac{\delta_k}{\eta2^{-l}}\right)^{\beta}\bigg)\int_{A(2\eta,\delta_k)}|\nabla\tilde{u}_k|^2\,,
 \end{align}
 where $\beta=|\log_2\mu|$.~\qed
\end{Prop}
\noindent The proof of Proposition \ref{propositionharmpart} is the same of that of of Lemma III.5 in \cite{DLGR}, therefore we omit it.

\noindent Next we recall the {\bf weighted Wente inequality} and the energy estimates that we have obtained \cite{DLGR}:
\begin{Lm}[Lemma F.1, \cite{DLGR}]\label{morrey1}
Let $a,b\in W^{1,2}(B_1)$ and $\varphi$ be a solution of
\begin{equation}
-\Delta\varphi=\p_{x_1}a\ \p_{x_2}b-\p_{x_2}a\ \p_{x_1}b\quad\mbox{in $B_1$}
\end{equation}
Then the following estimate holds.
\begin{eqnarray}\label{morrey2}
\int_{B_{\frac{1}{2}}\setminus B_{\frac{1}{4}}}|\nabla \varphi|^2\ dx&\le &\frac{2}{3}\int_{B_1\setminus B_{\frac{1}{2}}}|\nabla \varphi|^2+C\ \int_{B_1}|\nabla a|^2\ dx\ \int_{B_1} f(|x|)\ |\nabla b|^2\ dx
\end{eqnarray}
where $f$ is given by 
\be
 \label{weight-formula}
 f(r)=r^2\log^2\left(1+\frac{1}{r}\right)\log\left(1+\log\frac{1}{r}\right)\ ,
 \ee
 and $C>0$ is a universal constant. ~\qed
\end{Lm}

\begin{Lm}[Lemma F.2, \cite{DLGR}]\label{morrey9bis}
Under the hypothesis of lemma \ref{morrey1}, for every $\alpha\in (0,2)$ there exist $C_{\alpha}>0$ such that for any $k\in {\N}$
 \begin{eqnarray}\label{morrey9}
{\int_{{A}_{k+1}}}|\nabla \varphi|^2\ dx&\le&\gamma^{k+1}{\int_{ A_{0}}}|\nabla \varphi|^2\ dx+C_\al\ \int_{B_1}|\nabla a|^2\ dx\
\sum_{n=0}^{\infty}\gamma^{|n-k|}\int_{A_n}|\nabla b|^2\ dx\ .
\end{eqnarray}
where $\gamma=\max(2^{-\alpha},\frac{2}{3})$ and we are making use of the notation $A_k:=B_{2^{-k}}(0)\setminus B_{2^{-k-1}}(0)$.
\end{Lm}
\begin{Prop}\label{propositionwentepart}
Under the same assumptions of Lemma \ref{lemmahodgedec}, there exist $C>0$ and $\gamma\ge \frac{2}{3})$ such that
 \begin{equation*}
 \int_{A_j}|\nabla \varphi_k|^2\leq \gamma^{j+1}C\bigg(\int_{A(\eta,\delta_k)}|\nabla u_k|^2\bigg)^2+C_\alpha \bigg(\int_{A(\eta,\delta_k)}|\nabla u_k|^2\bigg)\sum\limits_{n=0}^\infty \gamma^{|n-j|}\int_{{A}_n}|\nabla \tilde{u}_k|^2.
 \end{equation*}
\end{Prop}
\noindent {\bf Proof of Proposition \ref{propositionwentepart}.} The proof of this proposition follows immediately from Lemmae \ref{morrey1} and \ref{morrey9bis} and is proved in full detail in Section III therein. 
\\ \noindent\qed\\ \noindent
\begin{Prop}\label{propabseries} Set 
\[a_j=\int_{\bar A_j}e_{\varepsilon_k}(u_k)\]
\begin{align}
 b_j=&C\int_{\tilde{A}_j}|\nabla h_k|^2+\,\gamma^{j+1}\tilde{C} \bigg(\int_{A(\eta,\delta_k)}|\nabla u_k|^2\bigg)^2 \\ \noindent&+\sum\limits_{|j'-j|\leq 3}C\sigma \int_{\bar{A}_{j'}}e_{\epsilon_k}(u_k)+\frac{C\varepsilon_k^2}{\sigma 2^{-2{j'}}}\int_{\bar{A}_{j'}}e_{\epsilon_k}(u_k)+\sigma \bigg(\frac{C2^{-2j}}{\epsilon_k^2}e^{-\frac{C2^{-j}}{\varepsilon_k}}\bigg)^2\,.
\end{align}
Then
\[a_k\leq b_k+\varepsilon_0\sum\limits_{n=0}^\infty\gamma^{|n-k|}a_n\,.\]
\end{Prop}
\noindent {\bf Proof of Proposition \ref{propabseries}.}
By equation \eqref{eee1} and Proposition \ref{propinneck},
\begin{align}
 \int_{\bar A_j}e_{\varepsilon_k}(u_k)&\leq C\int_{\tilde{A}_j}|\nabla h_k|^2+C\int_{\tilde{A}_j}|\nabla \varphi_k|^2 +\sum\limits_{|j'-j|\leq 3}C\sigma \int_{\bar{A}_{j'}}e_{\epsilon_k}(u_k)+\frac{C\varepsilon_k^2}{\sigma 2^{-2{j'}}}\int_{\bar{A}_{j'}}e_{\epsilon_k}(u_k)\\ \noindent+&\sigma \bigg(\frac{C2^{-2j}}{\epsilon_k^2}e^{-\frac{C2^{-j}}{\varepsilon_k}}\bigg)^2\,.
\end{align}
The claim then follows by Proposition \ref{propabseries}.
\\ \noindent\qed\\ \noindent

\noindent Combining Lemma \ref{lemmahodgedec}, \ref{lemma15}, \ref{lemmacontrolpotentialpart}, and Proposition \ref{propositionharmpart}, \ref{propositionwentepart} we can use Lemma G.1 in \cite{DLGR} and we obtain for every $\gamma\in (0,1)$ and $0<\mu<\gamma<1$, $2\delta_k/\eta\leq 2^{-s_2}\leq 2^{-s_1}\leq \eta/2$
\begin{equation}
 \sum\limits_{\ell=s_1}^{s_2}\mu^{|l-j|}a_{\ell}\leq \sum\limits_{\ell=s_1}^{s_2}\mu^{|l-j|}b_l+C_{\mu,\gamma}\varepsilon_0\sum\limits_{\ell=0}^{\infty}\mu^{|l-j|}a_{\ell}
\end{equation}
for some positive constants $C_{\mu,\gamma}$ independent of $a_k,b_k$.\\ \noindent
 This can be written as
 \begin{align}
 \sum\limits_{\ell=s_1}^{s_2}\mu^{|l-j|}\int_{\bar{A}_{\ell}}e_{\varepsilon_k}(u_k)\leq &C\bigg[\bigg(\frac{2^{-j}}{\eta}\bigg)^\beta+\bigg(\frac{\delta_k}{2^{-j}\eta}\bigg)^\beta\bigg]\int_{A(\eta,\delta_k)}|\nabla\tilde{u}_k|^2\\ \noindent
 +&C\bigg(\mu^j\bigg[\frac{\gamma}{\mu}\bigg]^{s_{1}}+\gamma^j\bigg)\bigg(\int_{A(\eta,\delta_k)}|\nabla\tilde{u}_k|^2\bigg)^2\\ \noindent+&\sigma C\sum\limits_{\ell=s_1}^{s_2}\mu^{|l-j|}\int_{\tilde{A}_l}e_{\varepsilon_k}(u_k)\\ \noindent+&\sum\limits_{l=s_1}^{s_2}\mu^{|l-j|}\sigma \bigg(\frac{C2^{-2l}}{\epsilon_k^2}e^{-\frac{C2^{-l}}{\varepsilon_k}}\bigg)^2\label{partwithsigma}\\ \noindent+&\sum\limits_{\ell=s_1}^{s_2}\mu^{|l-j|}\frac{C\varepsilon_k^2}{\sigma2^{-2l}}\int_{\bar{A}_l}e_{\varepsilon_k}(u_k) \\ \noindent
 +&\varepsilon_0C_{\mu,\gamma}\sum\limits_{\ell=0}^\infty\mu^{|l-j|}\int_{\bar{A}_j}e_{\varepsilon_k}(u_k)
\end{align}

\noindent The last two terms can be absorbed in the following way: notice that by Lemma \ref{lemmaepsilondeltaratio}, $\delta_k>>\varepsilon_k$ and for all $\ell\in \{s_1,...,s_2\}$, $2^{-2\ell}\geq \frac{\delta_k^2}{\eta^2}$, which implies
\begin{equation}
 \frac{\varepsilon_k^2}{\sigma 2^{-2l}}\leq\frac{\eta^2}{\sigma}\bigg(\frac{\varepsilon_k}{\delta_k}\bigg)^2\to 0 ~~\mbox{as $k\to +\infty$}\,.
\end{equation}
We proceed to bound the following term appearing in \eqref{partwithsigma}:
$$\sum\limits_{s_1\leq l\leq s_2}\mu^{|l-j|} \bigg(\frac{C2^{-2l}}{\epsilon_k^2}e^{-\frac{C2^{-l}}{\varepsilon_k}}\bigg)^2\,.$$
For $|x|\cong 2^{-j}$, using $\beta=|\log_2\mu|$, and $r^2e^{-Cr}\leq Cr^{-\beta/2}$ for all $r>>r_0$,
\begin{align}\label{extratermexp}
 &\sum\limits_{s_1\leq l\leq s_2}\mu^{|l-j|} \bigg(\frac{C2^{-2l}}{\epsilon_k^2}e^{-\frac{C2^{-l}}{\varepsilon_k}}\bigg)^2\\=& \frac{\delta_k^\beta}{|x|^\beta\eta^\beta}\sum\limits_{l=s_1}^{j}\mu^{|l-j|} \frac{|x|^\beta\eta^\beta}{\delta_k^\beta}\bigg(\frac{C2^{-2l}}{\epsilon_k^2}e^{-\frac{C2^{-l}}{\varepsilon_k}}\bigg)^2+\frac{|x|^\beta}{\eta^\beta}\sum\limits_{l=j}^{s_2}\mu^{|l-j|} \frac{\eta^\beta}{|x|^\beta}\bigg(\frac{C2^{-2l}}{\epsilon_k^2}e^{-\frac{C2^{-l}}{\varepsilon_k}}\bigg)^2\\ \noindent
 \leq& C\frac{\delta_k^\beta}{|x|^\beta\eta^{\beta}}\sum\limits_{l=s_1}^{j}\mu^{|l-j|}\eta^{2\beta} \frac{\epsilon_k^\beta}{\delta_k^\beta}2^{-\beta(j-l)}+\frac{|x|^\beta}{\eta^\beta}\sum\limits_{l=j}^{s_2}\mu^{|l-j|} \frac{\eta^\beta}{|x|^\beta}\epsilon_k^\beta2^{\beta l}\\ \noindent
 \leq& C\frac{\delta_k^\beta}{|x|^\beta\eta^{\beta}}\sum\limits_{l=s_1}^{j}\mu^{2|l-j|}\eta^{2\beta} \frac{\epsilon_k^\beta}{\delta_k^\beta}+\frac{|x|^\beta}{\eta^\beta}\sum\limits_{l=j}^{s_2}\mu^{l-j}\mu^j {\eta^{2\beta}}o_{\eta,k}(1)
 \\ \noindent\leq&Co_{\eta,k}(1)\bigg(\frac{\delta_k^\beta}{|x|^\beta\eta^\beta}+\frac{|x|^\beta}{\eta^\beta}\bigg)\,.
\end{align}
\noindent Inequality \eqref{extratermexp} allows us to absorb the terms in the estimate above to obtain
\begin{equation}\label{estekbeta}
 \int_{\bar{A}_{j}}e_{\varepsilon_k}(u_k)\leq \bigg[\bigg(\frac{2^{-j}}{\eta}\bigg)^\beta+\bigg(\frac{\delta_k}{2^{-j}\eta}\bigg)^\beta\bigg]o_{\eta,k}(1)
\end{equation}
for all $j\in\{s_1,...,s_2\}$.
\color{black}
 \section{$C^0$ estimates and $L^{2,1}$ quantization in the neck region}
 \begin{Th}\label{theoremneckcontrolc0}
Let $u_k:(\Sigma,h)\to \mathbb{R}^{n+1}$ be a sequence of critical points of ${E}_{\varepsilon_k}$ with uniformly bounded energy ${E}_{\varepsilon_k}(u_k)\leq \Lambda$ and assume that it converges in the bubble tree sense to $(u_\infty,v_\infty)$. Then there exist $C>0$, $0<\beta<1$ and $\eta_0>0$, such that for all $\eta<\eta_0$ and for all $x\in A(\eta,\delta_k):$
 \begin{equation}\label{pointwisecontrolintheneck}
 |x|^2|e_{\varepsilon_k}(u_k)|(x)\leq C\bigg(\frac{|x|^\beta}{\eta^\beta}+\frac{\delta_k^\beta}{\eta^\beta|x|^\beta}\bigg)o_{\eta,k}(1)
 \end{equation}
\end{Th}
\noindent {\bf Proof of Theorem \ref{theoremneckcontrolc0}.}
 Applying Proposition \ref{propereg} to \eqref{estekbeta} implies \eqref{pointwisecontrolintheneck}.
\\ \noindent\qed\\ \noindent

\noindent The pointwise estimate \eqref{pointwisecontrolintheneck} of the energy density permits us to get the following {\bf upgrade} of the energy quantization in the neck regions:
\begin{Co}\label{corollaryl21control}
 Let $u_k:(\Sigma,h)\to \mathbb{R}^{n+1}$ be a sequence of critical points of ${E}_{\varepsilon_k}$ with uniformly bounded energy ${E}_{\varepsilon_k}(u_k)\leq \Lambda$ and assume that it converges in the bubble tree sense to $(u_\infty,v_\infty)$. Then 
 \[\lim\limits_{\eta\to 0}\limsup\limits_{k\to\infty}||(e_{\epsilon_k}(u_k))^{1/2}||_{L^{2,1}(A(\eta,\delta_k))}=0\,.\]
\end{Co}
\noindent {\bf Proof of Corollary \ref{corollaryl21control}.}
 By \eqref{pointwisecontrolintheneck}, for $x\in A(\eta,\delta_k)$,
 \[|e_{\epsilon_k}(u_k)|\leq \frac{1}{|x|^2}\bigg(\frac{|x|^\beta}{\eta^\beta}+\frac{\delta_k^\beta}{\eta^\beta|x|^\beta}\bigg)o_{\eta,\delta_k}(1)\,.\]
 By definition of the $L^{2,1}$ norm,
 \begin{equation}
 ||(e_{\epsilon_k}(u_k))^{1/2}||_{L^{2,1}(A(\eta,\delta_k))}\approx \int_{0}^\infty\left(\mu\{x\in A(\eta,\delta_k)\,\,\vert\,\,(e_{\epsilon_k}(u_k))^{1/2}\geq t\}\right)^{1/2}\,.
 \end{equation}
 Now we observe that
 \begin{equation}
 |(e_{\epsilon_k}(u_k))^{1/2}|\le C \frac{1}{|x|}\bigg(\frac{|x|^{\beta/2}}{\eta^{\beta/2}}+\frac{(\delta_k)^{\beta/2}}{\eta^{\beta/2}|x|^{\beta/2}}\bigg)o_{\eta,\delta_k}(1).
 \end{equation}
 We have 
 \begin{align}
 &\left(\mu\{x\in A(\eta,\delta_k)\,\,\vert\,\,\frac{1}{|x|}\bigg(\frac{|x|^{\beta/2}}{\eta^{\beta/2}}+\frac{(\delta_k)^{\beta/2}}{\eta^{\beta/2}|x|^{\beta/2}}\bigg)\geq t\}\right)^{1/2}\nonumber\\ \noindent\leq&\left( \underbrace{\mu\{x\in A(\eta,\delta_k)\,\,\vert\,\, \frac{|x|^{\beta/2-1}}{\eta^{\beta/2}}\geq t\}}_{(1)}\right)^{1/2}+\left(\underbrace{\mu\{x\in A(\eta,\delta_k)\,\,\vert\,\,\frac{(\delta_k)^{\beta/2}}{\eta^{\beta/2}|x|^{1+\beta/2}}\geq t\}}_{(2)}\right)^{1/2}\,.
 \end{align}
{\bf First term $(1)$:}
 \begin{equation}
\mu\{x\in A(\eta,\delta_k)\,\,\vert\,\, \frac{|x|^{\beta/2-1}}{\eta^{\beta/2}}\geq t\}=\mu\{x\in A(\eta,\delta_k)\,\,\vert\,\,|x|\leq\left( \frac{1}{t\eta^{\beta/2}}\right)^{1/(1-\beta/2)}\}
 \end{equation}
 We can integrate with respect to $t$ between $ \frac{1}{\eta}$ and $\frac{\eta^{1-\beta}}{\delta_k^{1-\beta/2}}$:
 \begin{align}\label{firstterm}
 &\int_{ \frac{1}{\eta}}^{\frac{\eta^{1-\beta}}{\delta_k^{1-\beta/2}}}\left(\mu\{x\in A(\eta,\delta_k)\,\,\vert\,\,|x|\leq\left( \frac{1}{t\eta^{\beta/2}}\right)^{1/(1-\beta/2)}\}\right)^{1/2} dt \approx
\frac{1}{\eta^{\frac{\beta/2}{1-\beta/2}}}\bigg[t^{1-\frac{1}{1-\beta/2}}\bigg]_{ \frac{1}{\eta}}^{\frac{\eta^{1-\beta}}{\delta_k^{1-\beta/2}}}\leq C
 \end{align}
 
 \noindent
{\bf Second term $(2)$:} 
We can integrate with respect to $t$ between $ \frac{\delta_k^{\beta/2}}{\eta^{1+\beta}}$ and $\frac{\eta}{\delta_k}$:
 \begin{align}\label{secondterm}
 \mu\{x\in A(\eta,\delta_k)\,\,\vert\,\,\frac{\delta_k^{\beta/2}}{\eta^{\beta/2}|x|^{1+\beta/2}}\geq t\}\leq \mu\{x\in A(\eta,\delta_k)\,\,\vert\,\,|x|\leq\left(\frac{\delta_k^{\beta/2}}{t\eta^{\beta/2}}\right)^{1/(1+\beta/2)}\}
 \end{align}
 \begin{align}
 &\int_{\frac{\delta_k^{\beta/2}}{\eta^{1+\beta}}}^{\frac{\eta }{\delta_k }}\mu\{x\in A(\eta,\delta_k)\,\,\vert\,\,\frac{\delta_k^{\beta/2}}{\eta^{\beta/2}|x|^{1+\beta/2}}\geq t\}^{1/2}\leq C\left(\frac{\delta_k}{\eta}\right)^{\frac{\beta}{2+\beta}}\bigg[t^{1-\frac{1}{1+\beta/2}}\bigg]_{\frac{\delta_k^{\beta/2}}{\eta^{1+\beta}}}^{\frac{\eta}{\delta_k}} 
 \leq C\,. \end{align}
 
 This concludes the proof of the Corollary. \qed\\ \noindent

\section{Main theorem}
\noindent With the point-wise estimates obtained in Theorem \ref{theoremneckcontrolc0}, we can now follow the robust strategy developed in \cite{DLGR} to prove the upper-semi-continuity of the Morse index.\\ \noindent
\noindent The second variation of the Ginzburg-Landau functional is given by
 \begin{equation}\label{secondderivative}
 D^2 {E}_{\varepsilon_k}(u_k)[v]=\int_{\Sigma}|\nabla v|^2+\frac{2}{\varepsilon_k^2}\langle u_k,v\rangle^2-\frac{1}{\varepsilon_k^2}(1-|u_k|^2)|v|^2
 \end{equation}
 where $v\in W^{1,2}\cap\Gamma(u^{-1}TS^n)$. 
 By Corollary \ref{deltareglap}, 
 \begin{align}
 \|\Delta u_k\|_{L^\infty(A_j)}2^{-2j}&\leq C\int_{\tilde{A}_j}|\nabla u_k|^2+\frac{(1-|u_k|^2)^2}{\varepsilon_k^2}+C\frac{2^{-2j}}{\epsilon_k^2}e^{-C\frac{2^{-j}}{\epsilon_k}}
 \end{align}
and by \eqref{estekbeta},
 \begin{equation}
 \frac{1-|u_k(x)|^2}{\varepsilon_k^2}\leq C |x|^{-2}o_{\eta,k}(1)\bigg[\bigg(\frac{|x|}{\eta}\bigg)^{\beta}+\bigg(\frac{\delta_k}{\eta|x|}\bigg)^{\beta}\bigg]\,.\label{weightboundinannulus}
 \end{equation}
 Combining \eqref{secondderivative} and \eqref{weightboundinannulus},
 \begin{equation}
 D^2{E}_{\varepsilon_k}(u_k)[v]\geq Q_{u_k}(v)\ge \int_{A(\eta,\delta_k)}|\nabla v|^2-\frac{Co_{\eta,k}(1)}{|x|^2}\bigg[\bigg(\frac{|x|}{\eta}\bigg)^{\beta}+\bigg(\frac{\delta_k}{\eta|x|}\bigg)^{\beta}\bigg]|v|^2,
 \end{equation}
 for all $ v\in W^{1,2}_0(A(\eta,\delta_k))$ and
 \begin{equation}
 Q_{u_k}(v)=\int_\Sigma|\nabla v|^2-\frac{1}{\varepsilon_k^2}(1-|u_k|^2)|v|^2.
 \end{equation}
 We introduce the following weights on $\Sigma$ as in Section IV of \cite{DLGR},\vspace{0.3cm}
 \begin{equation}
 \omega_{\eta, k}=
 \begin{cases}
 \displaystyle{ \frac{1}{|x|^2}\bigg[\bigg(\frac{|x|}{\eta}\bigg)^{\beta}+\bigg(\frac{\delta_k}{\eta|x|}\bigg)^{\beta}\bigg] \text{ on }A(\eta,\delta_k)}\\[5mm]
 \displaystyle{ \frac{1}{\eta^2}\bigg[1+\bigg(\frac{\delta_k}{\eta^2}\bigg)^{\beta}\bigg] \text{ on }\Sigma\setminus B_\eta}\\[5mm]
 \displaystyle{ \frac{\eta^2}{\delta_k^2}\bigg[\frac{1}{\eta^4}\frac{(1+\eta^2)^2}{(1+\delta_k^{-2}|x|^2)^2}+\bigg(\frac{\delta_k}{\eta^2}\bigg)^{\beta}\bigg] \text{ on }B_{\delta_k/\eta}}
 \end{cases}
 \end{equation}
 By applying Lemma IV.1 in \cite{DLGR} we get the following result:
\begin{Lm}
\label{posinneck}
There exists $\la_0>0$ independent of $k$ and there exists $\eta_0$ such that for any $0<\eta<\eta_0$ and for $k$ large enough the following holds
\begin{equation*}
 \forall\ {w}\,\in W^{1,2}(\Sigma)\,\,\text{s.t.}\,\,{w}=0\,\, \text{in}\,\,\Sigma\setminus A(\eta,\delta_k)\,,\,\, D^2{E}_{\varepsilon_k}(u_k)[w]\geq Q_{u_k}(w)\geq \, \la_0\, \int_{\Sigma}\om_{\eta,k}\ w^2\ \,dvol_h\ .
\end{equation*}
\end{Lm}
 \noindent The proof of Lemma \ref{posinneck} is exactly the same as the one of Lemma IV.2 in \cite{DLGR}. We can now move state the main theorem.

\begin{Th}\label{morseindextheoremfordirichlet}
Let $u_k:(\Sigma,h)\to\mathbb{R}^{n+1}$ be a sequence of critical points of $ E_{\varepsilon_k}$ with uniformly bounded energy $ E_{\varepsilon_k}(u_k)\leq \Lambda$ and assume that it converges in the bubble tree sense to $(u_\infty,v_\infty)$. Then for $k$ large enough,
\begin{equation}
 \text{Ind}_{ E_{\varepsilon_k}}(u_k)+ \text{Null}_{ E_{\varepsilon_k}}(u_k)\leq \text{Ind}_{ E}(u_\infty)+ \text{Null}_{ E}(u_\infty)+\text{Ind}_{ E}(v_\infty)+\text{Null}_{ E}(v_\infty)\,.
\end{equation}
\end{Th}
\noindent {\bf Proof of Theorem \ref{morseindextheoremfordirichlet}.} We divide the proof of the Theorem into several steps.\vspace{0.3cm}\\ \noindent
\noindent \textbf{Step 1: Diagonalization with respect to the weights}\\ \noindent
The weighted $L^2$ scalar products with respect to the weights $ w_{\eta,k}$ are defined by
\[
<f,g>_{\om_{\eta,k}}=\int_\Sigma \ f\,g\ \om_{\eta,k}\, \,dvol_h\ .
\]
for all $f,g\in L^2_{\omega_{\eta, k}}:=\big\{w\in L^2(\Sigma),\, \int_\Sigma w^2 \omega_{\eta, k}\,dvol_h <+\infty\big\}$.
 
\noindent As in section IV.2 of \cite{DLGR} we are going to consider the diagonalization of suitable self-adjoint operators with respect to the weights $\om_{\eta,k}$:
\begin{eqnarray}\label{diag}
 Q_{u_k}(w)&=&\int_{\Sigma}|\nabla w|^2-\frac{1-|u_k|^2}{\varepsilon_k^2}|w|^2 \,dvol_h\\ \noindent
 &=&\langle \mathcal{L}_{\eta,k}w,w\rangle_{\omega_{\eta,k}}=\int_{\Sigma}\bigg[-\omega_{\eta,k}^{-1}\Delta_hw-\omega_{\eta,k}^{-1}\frac{1-|u_k|^2}{\varepsilon_k^2} w\bigg]w\omega_{\eta,k}\, \,dvol_h
\end{eqnarray}
where the weights are smooth and bounded $0<C_\eta<\omega_{\eta,k}< C_{\eta,k}$. By Lemma B.1 in \cite{DLGR} these operators are diagonalizable on $L^2_{\omega_{\eta, k}}\cap W^{1,2}$.

\begin{Lm}\label{indexdimrelation}
 Denoting $\mathcal{E}_{\eta,k}(\lambda)$ the eigenspace for the eigenvalue $\lambda$ for the operator $ \mathcal{L}_{\eta,k}$,
 \begin{equation}\label{firstestind}
 \mathrm{Ind}_{{E}_{\varepsilon_k}}(u_k)+ \mathrm{Nul}_{{E}_{\varepsilon_k}}(u_k)\leq\mathrm{dim}\bigg[\oplus_{\lambda\leq0}\mathcal{E}_{\eta,k}(\lambda)\bigg]=:\mathrm{dim}(\mathcal{E}^0_{\eta,k})
 \end{equation}
\end{Lm}
\noindent {\bf Proof of Lemma \ref{indexdimrelation}.}
 The estimate \eqref{firstestind} directly follows from the fact that by definition of $Q_{u_k}(w)$,
 \begin{equation}
 D^2{E}_{\varepsilon_k}(u_k)[w]\geq Q_{u_k}(w).
 \end{equation}
\\ \noindent\qed\\ \noindent
\noindent We now prove that the negative eigenvalues $\lambda$ for the operator $ \mathcal{L}_{\eta,k}$ are bounded from below, this is the analogue of Lemma IV.4 in the Ginzburg-Landau case:
\begin{Lm}
\label{lm-lower-spectr}
There exist two constants $\eta_0>0$ and $\mu_0>0$ independent of $k$ and a family of constants $\mu_{\eta,k}>0$ satisfying 
\be
\label{I.25-a}
\lim_{\eta\rightarrow 0}\limsup_{k\rightarrow +\infty} \mu_{\eta,k}=0\ \quad \mbox{ and }\quad \forall k\in{\N},~~\forall \eta \in(0,\eta_0)\quad\mbox{ one has }\quad 0<\mu_{\eta,k}<\mu_0\ ,
\ee
and such that for any element $\lambda\in{ \mathbb{R}}$ with
\be
\label{I.25-b}
\mbox{dim}\lf( {\mathcal E}_{\eta,k}(\la) \rg)>0\quad\Longrightarrow\quad\la\ge -\,\mu_{\eta,k}\ge - \mu_0\ .
\ee
\hfill $\Box$
\end{Lm}
\noindent {\bf Proof of Lemma \ref{lm-lower-spectr}.}
 By \eqref{weightboundinannulus},
 \begin{equation}
 \frac{1-|u_k|^2}{\varepsilon_k^2}\leq \omega_{\eta,k}(x)o_{\eta,k}(1),~~ \forall x\in A(\eta,\delta_k).
 \end{equation}
 We also have 
 \begin{equation}
 u_k\frac{1-|u_k|^2}{\varepsilon_k^2}=-\Delta u_k ~~\forall x\in\Sigma,
 \end{equation}
 which combined with $u_k\to u_\infty$ in $C^\infty_{\text{loc}}(\Sigma\setminus \{p\})$ and $v_k\to v_\infty$ implies
 \begin{equation}
 \bigg\vert\frac{1-|u_k|^2|}{\varepsilon_k^2}\bigg\vert\leq C\omega_{\eta,k}~~~ \forall x\in\Sigma.
 \end{equation}
 Define
 \begin{equation}
 \mu_{\eta,k}=\sup\limits_{x\in\Sigma}\bigg\vert\frac{(1-|{u}_k(x)|^2)}{\varepsilon_k^2\omega_{\eta,k}(x)}\bigg\vert
 \end{equation}
and let $\lambda$ be a negative eigenvalue for $ \mathcal{L}_{\eta,k}$. From the diagonalization of the operator we get $0\neq w$,
 \begin{equation}
\lambda = \lambda\int_\Sigma \omega_{\eta,k}w^2=\int_{\Sigma}|\nabla w|^2-\frac{1}{\varepsilon_k^2}(1-|u_k|^2)w^2\geq -\int_{\Sigma}\frac{1}{\varepsilon_k^2}(1-|u_k|^2)w^2\geq -\mu_{\eta,k}\int_\Sigma\omega_{\eta,k}w^2
 \end{equation}
 which concludes the proof of \eqref{I.25-b}. To prove \eqref{I.25-a} one proceeds exactly as in IV.28-IV.30 of \cite{DLGR}. This concludes the proof of the Lemma.
\\ \noindent\qed\\ \noindent
\noindent\textbf{Step 2: Limiting weights.} Define
\begin{align*}
 &\omega_{\eta,\infty}=\begin{cases}
 \frac{1}{\eta^2}\text{ on }\Sigma\setminus B_\eta\\ \noindent
 \frac{1}{|x|^2}\bigg(\frac{|x|}{\eta}\bigg)^{\beta} \text{ on }B_\eta
 \end{cases}\\ \noindent
 &\hat{\omega}_{\eta,\infty}=\begin{cases}
 \frac{1}{\eta^2}\bigg(\frac{1+\eta^2}{1+|y|^2}\bigg)^2\text{ on } B_{\eta^{-1}}\\ \noindent
 \frac{1}{\eta^{\beta}}\frac{1}{|y|^{2+\beta}} \text{ on }\mathbb{C}\setminus B_{\eta^{-1}}
 \end{cases}\\ \noindent
 &\tilde{\omega}_{\eta,\infty}(x)=\hat{\omega}_{\eta,\infty}(\pi(x))(1+|\pi(x)|^2)^2:S^2\to\mathbb{R}
\end{align*}
where $\pi:S^2\to\mathbb{C}$ is the stereographic projection.

\noindent The weights $\omega_{\eta,\infty}$ and $\hat{\omega}_{\eta,\infty}$ are defined on closed manifolds and satisfy the conditions of the diagonalization Lemma B.1 in \cite{DLGR}.\\ \noindent
We can diagonalize
\begin{equation}
 Q_{u_\infty}(w)=\langle \mathcal{L}_{\eta,\infty}w,w\rangle_{{\omega}_{\eta,\infty}}=\int_\Sigma [\omega_{\eta,\infty}^{-1}\Delta_hw-\omega_{\eta,\infty}^{-1}|\nabla u_\infty|^2w]w\omega_{\eta,\infty}
\end{equation}
on
\begin{equation}
 V_{u_\infty}=\bigg\{w\in L^2_{\omega_{\eta,\infty}}(\Sigma):\, \langle w,u_\infty\rangle = 0 \text{ a.e. on }\Sigma\bigg\}
\end{equation}
and denote 
\begin{equation}
 \mathcal{E}^0_{\eta,\infty}=\oplus_{\lambda\leq 0}\mathcal{E}_{\eta,\infty}(\lambda)\subset V_{u_\infty}.
\end{equation}
Since $D^2E(u_\infty)[w]=Q_{u_\infty}(w)$ for all $w\in L^2_{\omega_{\eta,\infty}}$,
\begin{equation}
 \text{dim}(\mathcal{E}^0_{\eta,\infty})\leq\text{Ind}_E(u_\infty)+\text{Nul}_E(u_\infty).
\end{equation}
We consider 
\begin{align*}
 &Q_{v_\infty}(w)=\langle \hat{ E}_{\eta,\infty}w,w\rangle_{\hat{\omega}_{\eta,\infty}}=\int_\mathbb{C} [\hat{\omega}_{\eta,\infty}^{-1}\Delta w-\hat{\omega}_{\eta,\infty}^{-1}|\nabla v_\infty|^2w]w\hat{\omega}_{\eta,\infty}=\\ \noindent
 &\tilde{Q}_{\tilde{v}_\infty}(\tilde{w})=\langle \tilde{ E}_{\eta,\infty}\tilde{w},\tilde{w}\rangle_{\tilde{\omega}_{\eta,\infty}}=\int_{S^2} [\tilde{\omega}_{\eta,\infty}^{-1}\Delta_{S^2} \tilde{w}-\tilde{\omega}_{\eta,\infty}^{-1}|\nabla \tilde{v}_\infty|^2\tilde{w}]\tilde{w}\tilde{\omega}_{\eta,\infty}.\label{117}
\end{align*}
By Lemma B.1 in \cite{DLGR}, $\tilde{Q}_{\tilde{v}_\infty}(\tilde{w})$ can be diagonalized on 
\begin{equation}
 V_{\tilde{v}_\infty}=\bigg\{w\in L^2_{\tilde{\omega}_{\eta,\infty}}(S^2):\, \langle w,\tilde{v}_\infty\rangle = 0 \text{ a.e. on }S^2\bigg\}.
\end{equation}
Let us prove the analogue of Proposition IV.1 in \cite{DLGR}.
\begin{Lm}\label{lemmam55}
 \begin{equation}
 \mathrm{dim}(\hat{\mathcal{E}}^0_{\eta,\infty})\leq\mathrm{Ind}_E(\tilde{v}_\infty)+\mathrm{Nul}_E(\tilde{v}_\infty).
 \end{equation}
\end{Lm}
\noindent {\bf Proof of Lemma \ref{lemmam55}.} See Proposition IV.1 in \cite{DLGR}\\ 
\qed
\\ \noindent
 So far we have defined diagonalizable operators on compact spaces satisfying
 \begin{align*}
 &\mathrm{Ind}_{E_{\varepsilon_k}}(u_k)+ \mathrm{Nul}_{E_{\varepsilon_k}}(u_k)\leq \mathrm{dim}(\mathcal{E}^0_{\eta,k}) \\ \noindent
 &\text{dim}(\mathcal{E}^0_{\eta,\infty})\leq\text{Ind}_E(u_\infty)+\text{Nul}_E(u_\infty)\\ \noindent
 & \mathrm{dim}(\hat{\mathcal{E}}^0_{\eta,\infty})\leq\mathrm{Ind}_E(\tilde{v}_\infty)+\mathrm{Nul}_E(\tilde{v}_\infty).
 \end{align*}
 To prove the Theorem it will therefore suffice to show the following Lemma.
 \begin{Lm}\label{lemmathatimpliesthm}
 \begin{equation}
 \mathrm{dim}(\mathcal{E}^0_{\eta,k})\leq \mathrm{dim}({\mathcal{E}}^0_{\eta,\infty})+\mathrm{dim}(\hat{\mathcal{E}}^0_{\eta,\infty})\text{ for }k \text{ large enough.}
 \end{equation}
 \end{Lm}
\noindent {\bf Proof of Lemma \ref{lemmathatimpliesthm}.}
 Consider the finite dimensional sphere
 \begin{equation}
 S_{\eta,k}=\bigg\{w\in\oplus_{\lambda\leq 0}\mathcal{E}_{\eta,k}(\lambda):\, \langle w,w\rangle_{\omega_{\eta,k}}=1\bigg\}
 \end{equation}
 Consider a sequence $w_k\in S_{\eta,k}$. Then since $\omega_{\eta,k}\geq \frac{1}{\eta^2}$ on $\Sigma$ for all $k$ big enough,
 \begin{equation}
 \int_{\Sigma}w_k^2\leq \eta^2.
 \end{equation}
 Moreover
 \begin{equation}
 0\geq \int_{\Sigma}|\nabla w_k|^2-\frac{1}{\varepsilon_k^2}(1-|u_k|^2)w_k^2\geq -\mu_0
 \end{equation}
 and
 \begin{equation}
 \int_{\Sigma}\frac{1}{\varepsilon_k^2}(1-|u_k|^2)w_k^2\leq \mu_0
 \end{equation}
 and we conclude that $w_k$ is bounded as a sequence in $W^{1,2}(\Sigma)$ and hence converges weakly
 \begin{equation}
 w_k\rightharpoonup w_\infty \in W^{1,2}(\Sigma,\mathbb{R}^{n+1}).
 \end{equation}
 Considering the blow ups
 \begin{equation}
 \hat{w}_k(y)=w_k(y\delta_k+x_k) \text{ on }\mathbb{C}
 \end{equation}
 \begin{equation}
 \int_{B_{\eta^{-1}}(0)}\hat{w}_k^2=\frac{1}{\delta_k^2}\int_{B_{\delta_k/\eta}}w_k^2\leq\frac{||\omega_{\eta,k}^{-1}||_{L^\infty(B_{\delta_k/\eta})}}{\delta_k^2}\leq C_\eta
 \end{equation}
 The gradients of $\hat{w}_k$ are bounded in $L^2$ by conformal invariance in two dimensions hence
 \begin{equation}
 \hat{w}_k\rightharpoonup \hat{w}_{\infty}\in W^{1,2}_{\text{loc}}(\mathbb{C},\mathbb{R}^{n+1}).
 \end{equation}
 \textbf{Step 3: Bootstrapping for stronger convergence.}\\ \noindent
 From the diagonalization properties discussed above, there are orthonormal families $\phi_k^j$ w.r.t. $L^2_{\omega_{\eta,k}}$ spanning $\mathcal{E}^0_{\eta,k}$ satisfying the eigenfunction equation for $ \mathcal{L}_{\eta,k}$
 \begin{equation}
 -\Delta \phi_k^j-\frac{1}{\varepsilon_k^2}(1-|u_k|^2)\phi_k^j=\lambda_k^j\phi_k^j\omega_{\eta,k},\quad \lambda_k^j\leq 0
 \end{equation}
 and for any element in the sphere $S_{\eta,k}$ there exist coefficients $c_k^j$ with
 \begin{equation}
 w_k=\sum\limits_{j=1}^{N_k}c_k^j\phi_k^j\quad \sum\limits_{j=1}^{N_k}(c_k^j)^2=1.
 \end{equation}
 We obtain the following equation for $w_k$
 \begin{equation}
 -\Delta w_k +\frac{1}{\varepsilon_k^2}(1-|u_k|^2)w_k=\sum\limits_{j=1}^{N_j}\lambda_k^jc_k^j\phi_k^j\omega_{\eta,k}
 \end{equation}
 and
 \begin{equation}
 \int_{\Sigma}\omega_{\eta,k}^{-1}|-\Delta w_k +\frac{1}{\varepsilon_k^2}(1-|u_k|^2)w_k|^2=\sum\limits_{j=1}^{N_j}(\lambda_k^j)^2(c_k^j)^2\leq (\mu_0)^2
 \end{equation}
 Moreover 
 \begin{equation}
 \sup\limits_k\bigg|\bigg|\frac{1}{\varepsilon_k^2}(1-|u_k|^2)\bigg|\bigg|_{L^\infty(\Sigma\setminus B_{\eta})}\leq C_\eta
 \end{equation}
 which by standard elliptic estimates implies up to subsequence
 \begin{equation}
 w_k\rightharpoonup w_\infty \text{ in }W^{2,2}_{\text{loc}}(\Sigma\setminus\{p\},\mathbb{R}^{n+1}).
 \end{equation}
 The exact same argument for the blow-ups also implies
 \begin{equation}
 \hat{w}_k\rightharpoonup\hat{w}_{\infty}\text{ in }W^{2,2}_{\text{loc}}(\mathbb{C}\setminus\{0\},\mathbb{R}^{n+1}).
 \end{equation}
 \textbf{Step 4: Nontrivial limit outside the necks.}\\ \noindent
 \textbf{Claim 1.} Either $w_\infty\neq 0$ or $\hat{w}_\infty\neq 0$.\\ \noindent
 \textbf{Proof of Claim 1.} We argue by contradiction. Suppose $w_\infty=0$ and $\hat{w}_\infty=0$. Then we can find a sequence $\check{w}_k\subset W^{1,2}_0(A(\eta,\delta_k),\mathbb{R}^{n+1})$ satisfying
 \begin{equation}
 \lim\limits_{k\to\infty}||\nabla w_k-\nabla \check{w}_k||_{L^2(\Sigma)}=0.
 \end{equation}
 Take for example the sequence
 \begin{equation}
 \check{w}_k=w_k\chi_k\bigg(2\frac{|x-x_k|}{\eta}\bigg)\bigg(1-\chi_k\bigg(\frac{\eta|x-x_k|}{\delta_k}\bigg)\bigg)
 \end{equation}
 where $\chi_k$ is a smooth approximation of the characteristic function $1_{B_1}$. Then
 \begin{align*}
 &|Q_{u_k}(w_k)-Q_{u_k}(\check{w}_k)|\leq \bigg\vert\int_{\Sigma\setminus B_{\eta/2}(x_k)}|\nabla w_k|^2-|\nabla\check{w}_k|^2-\frac{1}{\varepsilon_k^2}(1-|u_k|^2)(w_k^2-\check{w}_k^2)\bigg|+\\ \noindent&\bigg\vert\int_{B_{2\delta_k/\eta}(x_k)}|\nabla w_k|^2-|\nabla\check{w}_k|^2-\frac{1}{\varepsilon_k^2}(1-|u_k|^2)(w_k^2-\check{w}_k^2)\bigg|\to 0
 \end{align*}
 due to the assumption $w_\infty=0, \hat{w}_\infty=0$ as well as $\text{supp}(\check{w}_k)\subset A(\eta,\delta_k)$. Moreover
 \begin{equation}
 \lim\limits_{k\to\infty}\int_\Sigma \omega_{\eta,k}\check{w}_k^2=1.
 \end{equation}
 By Lemma \ref{posinneck}, 
 \begin{equation}
 Q_{u_k}(\check{w}_k)\geq c,
 \end{equation}
 but this contradicts $Q_{u_k}(w_k)\leq 0$ and therefore concludes the proof of the claim.\\ \\\noindent
 \textbf{End of the proof of Theorem \ref {morseindextheoremfordirichlet}.}\\ \noindent
 Let $N=\limsup\limits_k\mathrm{dim}(\mathcal{E}^0_{\eta,k})$ and $(\phi_k^j)_{j=1,...,N}$ families of negative eigenfunctions as before. By the strong convergence proved above we can assume,
 \begin{equation}
 \phi_k^j\rightharpoonup \phi_\infty^j \text{ in } W^{2,2}_{\text{loc}}(\Sigma\setminus\{p\})
 \end{equation}
 \begin{equation}
 \hat{\phi}_k^j\rightharpoonup \hat{\phi}_\infty^j \text{ in } W^{2,2}_{\text{loc}}(\mathbb{C}\setminus\{p\})
 \end{equation}
 Notice that by weak convergence of $u_k\to u_\infty$ in $W^{1,2}(\Sigma)$ and strong convergence of the blow ups $\hat{u}_k\to\hat{u}_\infty$ in $C^\infty_{\text{loc}}(\mathbb{C}),$
 \begin{equation}
 \frac{1}{\varepsilon_k^2}(1-|u_k|^2)\to |\nabla u_\infty|^2\text{ in }\mathcal{D}'
 \end{equation}
 \begin{equation}
 \frac{\delta_k^2}{\varepsilon_k^2}(1-|\hat{u}_k|^2)\to |\nabla \hat{u}_\infty|^2\text{ in }\mathcal{D}'\,.
 \end{equation}
 We therefore obtain equations solved in the weak sense by our limits of basis vectors,
 \begin{equation}
 \Delta \phi_\infty^j-|\nabla u_\infty|^2\phi_\infty^j=\lambda_\infty^j\omega_{\eta,\infty}\phi_\infty^j
 \end{equation}
 \begin{equation}
 \Delta \hat{\phi}_\infty^j-|\nabla v_\infty|^2\phi_\infty^j=\lambda_\infty^j\hat{\omega}_{\eta,\infty}\hat{\phi}_\infty^j
 \end{equation}
 for some $\lambda_\infty^j\leq 0$ where we use the previous notation $\hat{u}_\infty=v_\infty$. By Claim 1, for each $j=1,...,N$, 
 \begin{equation}
 (\phi_\infty^j,\hat{\phi}_\infty^j)\neq (0,0).
 \end{equation}
 Assume by contradiction
 \begin{equation}
 N>\mathrm{dim}({\mathcal{E}}^0_{\eta,\infty})+\mathrm{dim}(\hat{\mathcal{E}}^0_{\eta,\infty}).
 \end{equation}
 Then $(\phi_\infty^j,\hat{\phi}_\infty^j)$ must be linearly dependent, i.e. there is a non-zero family of coefficients s.t.
 \begin{equation}
 \sum\limits_{j=1}^Nc_\infty^j\phi_\infty^j=0\quad \text{and} \quad \sum\limits_{j=1}^Nc_\infty^j\hat{\phi}_\infty^j=0.
 \end{equation}
 Considering $w_k:=\sum\limits_{j=1}^Nc_\infty^j\phi_k^j$ we would have a sequence converging to zero both at the macroscopic ($w_k$) and the microscopic $\check{w}_k$ level, which contradicts Claim 1. This concludes the proof of Lemma \ref{lemmathatimpliesthm} and therefore of Theorem \ref{morseindextheoremfordirichlet}.\\\qed
 
 \appendix
 \section{proof of small-energy regularity}
 \noindent In this appendix we prove the small-energy regularity property for the Ginzburg-Landau perturbation, that is Proposition \ref{propereg}.

An important underlying tool is 
 {\em Bochner's identity} that we first present in the following lemma.
 
 \begin{Lm}[Bochner's identity]\label{Bochner}
 Let $u\in W^{1,2}(\Sigma,\mathbb{R}^{n+1})$ be a critical point of $E_\varepsilon$. Then the energy density $e_\varepsilon(u)$ satisfies
 
 \begin{eqnarray}
 - \Delta_h(e_\varepsilon(u)) 
 &\le & C (e_\varepsilon^2(u)+ e_\varepsilon(u))~~\mbox{in $\Sigma$}
 \end{eqnarray}
 for some positive constant $C$ which is independent of $u$ and $\varepsilon$.
\end{Lm}
\noindent {\bf Proof of Lemma \ref{Bochner}.}

To simplify the computation we take isothermal coordinates so that we can assume that $h$ is conformal to the flat metric, namely
$h_{ij}=e^{2\lambda}\delta_{ij}$ for some smooth function $\lambda.$\par
We set 
\begin{equation}\label{density}
e_{\varepsilon}(u)=\frac{1}{2}|d u|_h^2+\frac{1}{4\varepsilon^2}(1-|u|^2)^2
\end{equation}
We recall that
\begin{equation}\label{lapdensity}\Delta_h e_{\varepsilon}(u)=-e^{-2\lambda}\Delta e_{\varepsilon}(u)
\end{equation}
and
\begin{equation}\label{graddensity}
|d u|_h^2=e^{-4\lambda}|\nabla u|^2.\end{equation}
{\bf 1. } We have
\begin{eqnarray}\label{Bochner1}
\frac{1}{2}\Delta(e^{-4\lambda}|\nabla u|^2)&=&\frac{1}{2}(\Delta e^{-4\lambda})|\nabla u|^2+\nabla e^{-4\lambda}\cdot \nabla |\nabla u|^2\nonumber\\ \noindent
&=&\frac{1}{2}e^{-4\lambda}\Delta |\nabla u|^2.
\end{eqnarray}
Now we have
\begin{eqnarray}\label{Bochner2}
\frac{1}{2} \Delta |\nabla u|^2&=&\partial_{ii}^2\big(\frac{1}{2}\partial_ku\partial_ku\big)=\partial_i (\partial_{ik}^2u\partial_k u)\nonumber\\ \noindent
&=&|\nabla^2u|^2+\nabla\Delta u\cdot\nabla u.
\end{eqnarray}
Using the fact that $u$ is a critical point and therefore solves the Euler-Lagrange equation
\begin{equation}
 - \Delta u = u\ \frac{1}{\varepsilon^2}(1-|u|^2)
\end{equation}
we can compute $\nabla\Delta u\cdot\nabla u$ and we obtain
\begin{eqnarray}
 - \nabla\Delta u\cdot\nabla u&=&\nabla\big( \frac{1}{\varepsilon^2}(1-|u|^2)u\big)\cdot\nabla u=\nonumber\\ \noindent
 &=& \frac{1}{\varepsilon^2}\partial_i (1-|u|^2)u^j\partial_iu^j=\nonumber\\ \noindent
 &=& \frac{1}{\varepsilon^2}(-2u^k\partial_iu^ku^j\partial_iu^j)+(1-|u|^2)\partial_iu^j\partial_iu^j)=\nonumber\\ \noindent
 &=&-\frac{2}{\varepsilon^2}(u^j\partial_iu^j)^2+\frac{1}{\varepsilon^2}(1-|u|^2)|\nabla u|^2.\label{Bochner3}
\end{eqnarray}
{\bf 2.} We compute
$\Delta(\displaystyle\frac{1}{4\varepsilon^2}(1-|u|^2)^2)$.\par
\begin{eqnarray}\label{Bochner4}
\Delta(\frac{1}{4\varepsilon^2}(1-|u|^2)^2)&=& \frac{1}{4\varepsilon^2}\partial_{ii}^2\big((1-u^ju^j)^2\big)\nonumber\\ \noindent
&=&\frac{2}{\varepsilon^2}(-u^j\partial_iu^j)^2-\frac{1}{\varepsilon^2}(1-|u|^2)(\partial_iu^j\partial_iu^j+u^j\partial^2_{ii}u^j)\nonumber\\ \noindent
&=&\frac{2}{\varepsilon^2}(u^j\partial_iu^j)^2-\frac{1}{\varepsilon^2}(1-|u|^2)(|\nabla u|^2+u\cdot\Delta u)\nonumber\\ \noindent
&=&\frac{2}{\varepsilon^2}(u^j\partial_iu^j)^2-\frac{1}{\varepsilon^2}(1-|u|^2)(|\nabla u|^2-\frac {1}{\varepsilon^2}(1-|u|^2)|u|^2),
\end{eqnarray}
where in the last step in \eqref{Bochner4} we use the fact that $u$ solves the Euler-Lagrange equation:
\begin{equation}
 u\cdot\Delta u= -\frac {1}{\varepsilon^2}(1-|u|^2)|u|^2.\label{qy}
\end{equation}

The combination of the estimates 
\eqref{Bochner1}-\eqref{Bochner4} yields
\begin{eqnarray}
 - \Delta(e_\varepsilon(u))&=&-e^{-4\lambda}|\nabla^2 u|^2-e^{-4\lambda}\frac{4}{\varepsilon^2}(u^j\partial_iu^j)^2+e^{-4\lambda}\frac{2}{\varepsilon^2}(1-|u|^2)|\nabla u|^2-
 e^{-4\lambda} \frac{1}{\varepsilon^4}(1-|u|^2)^2|u|^2\nonumber\\ \noindent 
 &+&\frac{1}{2}(\Delta e^{-4\lambda})|\nabla u|^2+\nabla e^{-4\lambda}\cdot \nabla |\nabla u|^2\nonumber\\ \noindent
 &\le& -e^{-4\lambda}|\nabla^2 u|^2+e^{-4\lambda}\frac{2}{\varepsilon^2}(1-|u|^2)|\nabla u|^2+C_\lambda\nabla |\nabla u|^2+C_\lambda|\nabla u|^2.\label{Bochner5}
\end{eqnarray}
 To conclude the proof we show that it is possible to bound (point-wise in $x$) the last term in \eqref{Bochner5} in terms of the energy density $ e_\varepsilon(u)$ and of $e_\varepsilon(u)^2$.
 
\begin{eqnarray}
 \frac{2}{\varepsilon^2}(1-|u|^2)|\nabla u|^2-|\nabla ^2 u|^2&\leq& \frac{2}{\varepsilon^2}(1-|u|^2)|u||\nabla u|^2+\frac{2}{\varepsilon^2}(1-|u|^2)(1-|u|)|\nabla u|^2-|\nabla ^2 u|^2\nonumber\\ \noindent
 & {\leq}&2|\Delta u||\nabla u|^2+\frac{2}{\varepsilon^2}(1-|u|^2)^2|\nabla u|^2-|\nabla ^2 u|^2\nonumber\\ \noindent
 &{\leq} &\delta_1|\Delta u |^2 + 2C_{\delta_1}|\nabla u|^4 +\bigg(\frac{1}{\varepsilon^2}\big(1-|u|^2\big)^2\bigg)^2-|\nabla^2u|^2\nonumber\\ \noindent
 & \leq& Ce_\varepsilon(u)^2,\label{Bochner6}
\end{eqnarray}
where we used in \eqref{Bochner6} that $u$ solves the Euler Lagrange equation and $(1-|u|)\leq (1-|u|^2)$ by Lemma \ref{supestimatelemma}.\par
We also have
\begin{equation}
\nabla |\nabla u|^2= 2\nabla ^2 u\cdot \nabla u\le \delta_2|\nabla ^2 u|^2+C_{\delta_2}|\nabla u|^2.
\end{equation}
Since $|\Delta u|^2\leq |\nabla ^2 u|^2$, we can choose $\delta_1$ and $\delta _2$ in such a way that 
$$e^{-4\lambda}\delta_1|\Delta u |^2-e^{-4\lambda}|\nabla^2 u|^2+C_\lambda\delta_2|\nabla ^2 u|^2<0.$$
Hence from \eqref{lapdensity}, \eqref{Bochner5} and \eqref{Bochner6} it follows 
\begin{eqnarray}
 - \Delta_h(e_\varepsilon(u))&\le& C e^{-2\lambda} e_\varepsilon^2(u)+C|\nabla u|^2\nonumber\\ \noindent
 &\le & C_\lambda (e_\varepsilon^2(u)+ e_\varepsilon(u))
 \end{eqnarray}
where $C_\lambda$ is a positive constant depending on $\lambda, \nabla \lambda $ and $D^2\lambda.$
 This concludes the proof of the lemma.~~~$\Box$

We now prove Proposition \ref{propereg}, which thanks to Lemma \ref{Bochner} follows the lines of the proof for the harmonic map case in \cite{Sch} and we reproduce here for the reader's convenience.

\noindent {\bf Proof of Proposition \ref{propereg}.}
Let $u$ be a critical point of $E_\varepsilon$ and $r>0$ be such that
\begin{equation}\label{delta0i}
 \int_{B_r}\frac{1}{2}|\nabla u|^2+\frac{(1-|u|^2)^2}{4\varepsilon^2}<\delta_0
\end{equation}
where $\delta_0$ will be fixed later.
 Since the critical points of $E_\varepsilon$ are smooth,
\begin{align*}
 \sigma\in(0,\rho)\mapsto (r-\sigma)^2\max\limits_{x\in \overline{B}_\sigma}[{e}_\varepsilon(u)(x)]
\end{align*}
is continuous and we can define $\sigma_0$ such that
\begin{align*}
 (r-\sigma_0)^2\max\limits_{x\in \overline{B}_{\sigma_0}}[{e}_\varepsilon(u)(x)]=\max\limits_{\sigma\in [0,r]}(r-\sigma)^2\max\limits_{x\in \overline{B}_\sigma}[{e}_\varepsilon(u)(x)].
\end{align*}
By continuity of ${e}_\varepsilon(u)$ we can define
\begin{equation}
e_0=\max\limits_{x\in\overline{B}_{\sigma_0}}[{e}_\varepsilon(u)(x)]={e}_\varepsilon(u)(x_0)
\end{equation}
Setting $\rho_0=\frac{1}{2}(r-\sigma_0)$,
\begin{align*}
 \max\limits_{\overline{B}_{\rho_0}(x_0)}[{e}_\varepsilon(u)(x)]\leq \max\limits_{\overline{B}_{\sigma_0+\rho_0}(0)}[{e}_\varepsilon(u)(x)]\leq \bigg(\frac{r-\sigma_0}{r-(\rho_0+\sigma_0)}\bigg)^2 e_0\leq 4e_0.
\end{align*}
Define the re-scaled function
\begin{equation}
 v(y)=u\bigg(\frac{y}{\sqrt{e_0}}+x_0\bigg)
\end{equation}
and set $r_0=\sqrt{e_0}\rho_0$. Then we have
\begin{eqnarray}\label{RegP1}
 \sup\limits_{y\in B_{r_0}(0)}{e}_\varepsilon(v)(y)&
 = & \sup\limits_{y\in B_{r_0}(0)}\frac{1}{2e_0} |\nabla u|^2\bigg(\frac{y}{\sqrt{e_0}}+x_0\bigg)+\frac{1}{4\varepsilon^2}\bigg(1-\bigg\vert u\bigg(\frac{y}{\sqrt{e_0}}+x_0\bigg)\bigg\vert^2\bigg)^2\nonumber
 \\ \noindent&=&\frac{1}{e_0}\sup\limits_{y\in B_{r_0}(0)}{e}_{\varepsilon/\sqrt{e_0}}u\bigg(\frac{y}{\sqrt{e_0}}+x_0\bigg) 
\end{eqnarray}
and hence setting $\eta:=\varepsilon\sqrt{e_0}$ we have
\begin{equation}
 {e}_\eta(v)(0)= 1
\end{equation}
\begin{equation}
 \sup\limits_{y\in B_{r_0}(0)}{e}_\eta(v)(y)\leq 4.
\end{equation}
 \underline{Claim 1:} $r_0 < 1.$ \par
\noindent\textit{Proof of the Claim 1:} 
 If we assume by contradiction $r_0\geq 1$. Then $\sqrt{e_0}\ge 1$ as well.
 From Lemma \ref{Bochner} it follows that 
\begin{eqnarray}\label{equv}
\Delta_y {e}_\varepsilon(v)(y)&=&\frac{1}{e^2_0}[\Delta_x {e}_{\varepsilon/\sqrt{e_0}}(u)]\bigg(\frac{y}{\sqrt{e_0}}+x_0\bigg)\nonumber\\ \noindent
&\ge &-C\frac{1}{e^2_0}\left({e}_{\varepsilon/\sqrt{e_0}}(u)^2\bigg(\frac{y}{\sqrt{e_0}}+x_0\bigg)+{e}_{\varepsilon/\sqrt{e_0}}(u)\bigg(\frac{y}{\sqrt{e_0}}+x_0\bigg)\right)\nonumber
\\ \noindent
&\ge& -C\frac{1}{e^2_0} {e}_{\varepsilon/\sqrt{e_0}}(u)\bigg(\frac{y}{\sqrt{e_0}}+x_0\bigg)-C\frac{1}{e^2_0}{e}_{\varepsilon/\sqrt{e_0}}(u)\bigg(\frac{y}{\sqrt{e_0}}+x_0\bigg)\nonumber\\ \noindent
&\ge&-C {e}_\varepsilon(v)(y) .\end{eqnarray}

Since ${e}_\varepsilon(v)(y)$ is a subsolution of the equation $-\Delta w-C w=0$, by applying Harnack inequality to ${e}_\varepsilon(v)(y)$ in $B_1$
(by observing that $B_1\subset B_{r_0}$)
 we have
\begin{equation}\label{estv1}
 {e}_\eta(v)(0)= 1\leq C\int_{B_1}{e}_\eta(v)
\end{equation}
and by scaling,
\begin{eqnarray}\label{estv2}
 \int_{B_1(0)}{e}_\eta(v)&=& \int_{B_{\frac{1}{\sqrt{e_0}}}(x_0)}e_\varepsilon (u)
 \le \int_{B_{\rho_0}(x_0)}e_\varepsilon (u)\le \int_{B_{\sigma_0+\rho_0}(0)}e_\varepsilon (u)\nonumber\\ \noindent
 &\le & \int_{B_{r}(0)}e_\varepsilon (u) \leq \delta_0,
\end{eqnarray}
where we used that $r_0=\sqrt{e_0}\rho_0\geq 1$ implies $\sqrt{e_0}^{-1}\leq\rho_0\leq r$, $\sigma_0+\rho_0\le r$ and the hypothesis
\begin{equation}
 \int_{B_r}e_\varepsilon(u)<\delta_0\,.
\end{equation}
If $\delta_0$ is chosen small enough by combining \eqref{estv1} and \eqref{estv2} we get a contradiction and this proves the Claim 1.
\par
Since $r_0\le 1$ then $\sqrt{e_0}<\rho_0^{-1}$. We set $w(y)=u(x_0+\rho_0 y)$ and using Lemma \ref{Bochner} we get:
\begin{eqnarray}\label{equv2}
\Delta_y {e}_\varepsilon(w)(y)&=&\rho_0^2 [\Delta_x {e}_{\rho_0\varepsilon}(u)] (\rho_0 y+x_0 )\nonumber\\ \noindent
&\ge& -C\rho_0^4 [e^2_{\rho_0\varepsilon} (u) (\rho y+x_0 )+e_{\rho_0\varepsilon} (u) (\rho_0 y+x_0 )]\nonumber\\ \noindent
&\ge& -C\rho_0^2e_{\rho_0\varepsilon} (u)(\rho_0 y+x_0)=-Ce_\varepsilon (w).
\end{eqnarray}
where in \eqref{equv2} we used the fact that $\rho_0^2e_{\rho_0\varepsilon} (u) (\rho_0 y+x_0 )\le 1$.\par

 By applying Harnack to $w$ we get
\begin{eqnarray}\label{estr2}
 {e}_{\varepsilon}(w)(0)&=&r^2_0= \rho_0^2e_0\le C \int_{B_{1}(0)}{e}_{\varepsilon}(w)(y)=\rho^2_0\int_{B_{1}(0)}{e}_{\rho_0\varepsilon}(u)(\rho_0 y+x_0 ) \nonumber\\ \noindent
 &=&C \int_{B_{\rho_0}(x_0)}{e}_{\rho_0\varepsilon}(u)(x)\le C \int_{B_{r}(0)}{e}_{\rho_0\varepsilon}(u)(x).
\end{eqnarray}
Hence we get:
\begin{equation}\label{estr2bis} 
 \rho_0^2e_0\le C \int_{B_{r}(0)}{e}_{\rho_0\varepsilon}(u)(x).\end{equation}
 By the choice of $\sigma_0$, the estimate \eqref{estr3} implies
 \begin{eqnarray}
 (r-\sigma_0)^2\sup\limits_{\bar{B}_{\sigma_0}(0)}e_\varepsilon(u)&\leq &(r-\sigma_0)^2e_0\le 4 \rho_0^2e_0\le 4C \int_{B_r(0)}{e}_\varepsilon(u)
\end{eqnarray}
 and therefore 
\begin{equation}\label{estr3}
 \max\limits_{\sigma}(r-\sigma)^2\sup\limits_{B_{\sigma}(0)}e_\varepsilon(u)\leq 4C \int_{B_r(0)}{e}_\varepsilon(u)
\end{equation}
By choosing in \eqref{estr3} $\sigma=\frac{1}{2}r$ we obtain the estimate \eqref{Linfty} and we can conclude the proof of Proposition \ref{propereg}.~~~$\Box$\par
\bigskip

\begin{Rm}\label{Goodp}
Let $u_k$ be a sequence of critical points to $E_{\varepsilon_k}$ and $x\in\Sigma$ be such that there exists $\rho>0$ satisfying
\begin{equation}\label{goodp}
 \int_{B_\rho(x)}e_{\varepsilon_k}(u_k)\,dx\leq \delta_0
\end{equation}
for all $k$. Proposition \ref{propereg} implies that 
\begin{equation}
 ||\nabla u_k||_{L^\infty(B_{\rho/2}(x))}^2\leq C\frac{1}{{\rho^2}}{\int_{B_\rho(x)}e_{\varepsilon_k}(u_k)\,dx} 
\end{equation}
 and in particular $u_k$ is uniformly bounded in $W^{1,p}(B_{\rho/2}(x))$ for any $p\leq\infty$ and there is a sub-sequence converging strongly in $C^{\eta}(B_{\rho/2}(x))$ to a map $u_\infty:B_{\rho/2}(x)\to\mathbb{R}^{n+1}$. Moreover, due to weak convergence,
 \begin{equation}
 u_k\wedge \Delta u_k=0 \implies u_\infty\wedge \Delta u_\infty=0\,,
 \end{equation}
 which together with $|u_\infty|\equiv 1$ implies that $u_\infty$ is a harmonic maps taking values in $S^{n-1}$. In particular it must hold
 \begin{equation}
 -\Delta u_\infty = u_\infty |\nabla u_\infty|^2\quad in\,\, B_{\rho/2}(x)\,.~~~\Box
 \end{equation}
 \end{Rm}


\begin{thebibliography}{99}
 \bibitem{ChMa} Chodosh, Otis; Mantoulidis, Christos {\em Minimal surfaces and the Allen-Cahn equation on 3-manifolds: index, multiplicity, and curvature estimates.} Ann. of Math. (2) 191 (2020), no. 1, 213-328. 
 \bibitem{DLGR} Da Lio Francesca, Gianocca Matilde., Rivi\`ere Tristan, {\em Morse Index Stability for Critical Points to Conformally invariant Lagrangians}, preprint arXiv:2212.03124v2
 \bibitem{DLRp} Da Lio, Francesca; Rivi\`ere, Tristan {\em Conservation Laws for $p$-Harmonic Systems with Antisymmetric Potentials and Applications.} arXiv: 2311.04029.
\bibitem{dPKWY} Del Pino, Manuel; Kowalczyk, Michal ; Wei, Juncheng ; Yang, Jun 
{\em Interface foliation near minimal submanifolds in Riemannian manifolds with positive Ricci curvature. }
Geom. Funct. Anal. 20 (2010), no. 4, 918-957.
\bibitem{Hel} H\'elein, Fr\'ed\'eric {\sc Harmonic maps, conservation laws and moving frames. }Translated from the 1996 French original. With a foreword by James Eells. Second edition. Cambridge Tracts in Mathematics, 150. Cambridge University Press, Cambridge, 2002. 
 \bibitem{KaSt2} Karpukhin, M., Stern, D. Existence of harmonic maps and eigenvalue optimization in higher dimensions. Invent. math. 236, 713–778 (2024). https://doi.org/10.1007/s00222-024-01247-3
 \bibitem{KaSt} Karpukhin Mikhail, Stern Daniel, {\em Min-max harmonic maps and a new characterization of conformal eigenvalues.} arXiv:2004.04086
 \bibitem{Lam} Lamm, Tobias. “Energy identity for approximations of harmonic maps from surfaces.” Transactions of the American Mathematical Society, vol. 362, no. 8, 2010, pp. 4077–97. JSTOR, http://www.jstor.org/stable/25733357
 \bibitem{LaRi} Laurain, Paul; Rivi\`ere, Tristan {\em Angular energy quantization for linear elliptic systems with antisymmetric potentials and applications.} Anal. PDE 7 (2014), no. 1, 1-41.
 \bibitem{LiRi1} Lin, Fang-Hua; Rivi\`ere, Tristan {\em A quantization property for static Ginzburg-Landau vortices.} Comm. Pure Appl. Math., 54: 206-228. 
\bibitem{LiWa} Li, Yuxiang, Wang, Youde. (2015). A counterexample to the energy identity for sequences of approximate-harmonic maps. Pacific Journal of Mathematics. 274. 107-123. 10.2140/pjm.2015.274.107. 
\bibitem{LinWa} Lin, Fang-Hua; Wang, Changyou {\em Harmonic and Quasi-Harmonic Spheres}, Communications in Analysis and Geometry, Volume 7, Number 2, 397-429, 1999.
\bibitem{LinWa2} Lin, Fang-Hua; Wang, Changyou {\em Harmonic and Quasi-Harmonic Spheres, Part II}, Communications in Analysis and Geometry, Volume 10, Number 2, 341-375, 2002.
 \bibitem{LiZh} Li, Jiayu; Zhu, Xiangrong, {\em
Energy identity and necklessness for a sequence of Sacks-Uhlenbeck maps to a sphere}, Ann. Inst. H. Poincar\'e C Anal. Non Lin\'eaire (2019), no.1, 103-118.
\bibitem{MaNe} Marques, Fernando C.; Neves, Andr\'e {\em Morse index of multiplicity one min-max minimal hypersurfaces.} Adv. Math. 378 (2021), Paper No. 107527, 58 pp.
\bibitem{Mi} Michelat, Alexis; {\em Morse Index Stability of Biharmonic Maps in Critical Dimension}. https://arxiv.org/abs/2312.07494.
\bibitem{MiRi} Michelat, Alexis; Rivi\`ere, Tristan {\em Pointwise Expansion of Degenerating Immersions of Finite Total Curvature}. J. Geom. Anal. 33 (2023), no. 1, 24.
\bibitem{MiRi2} Michelat, Alexis; Rivi\`ere, Tristan {\em Morse Index Stability of Willmore Immersions I}. https://arxiv.org/abs/2306.04608 .
\bibitem{MiRi3} Michelat, Alexis; Rivi\`ere, Tristan {\em Weighted Eigenvalue Problems for Fourth-Order Operators in Degenerating Annuli}, https://arxiv.org/abs/2306.04609 \bibitem{Moo} Moore, John Douglas {\em Introduction to global analysis. Minimal surfaces in Riemannian manifolds.} Graduate Studies in Mathematics, 187. American Mathematical Society, Providence, RI, 2017. 
 \bibitem{MoRe} Moore, John Douglas and Ream, Robert {\em Minimal two-spheres of low index in manifolds with positive complex sectional curvature. }
Math. Z. 291 (2019), no. 3-4, 1295-1335. 
 \bibitem{Riv2} Rivi\`ere, Tristan {\em Lower semi-continuity of the index in the viscosity method for minimal surfaces.} Int. Math. Res. Not. IMRN 2021, no. 8, 5651-5675.
 \bibitem{SaUh} Sacks, J., and K. Uhlenbeck. “The Existence of Minimal Immersions of 2-Spheres.” Annals of Mathematics, vol. 113, no. 1, 1981, pp. 1–24. JSTOR, https://doi.org/10.2307/1971131. 
 \bibitem{Sch} Schoen, Richard, {\em Analytic Aspects of the Harmonic Map Problem.} In: Chern, S.S. (eds) Seminar on Nonlinear Partial Differential Equations. Mathematical Sciences Research Institute Publications, vol 2. Springer, New York, NY. 
\bibitem{Wen} Wente, H. "A General Existence Theorem for Surfaces of Constant Mean Curvature.." Mathematische Zeitschrift 120 (1971): 277-288.
\bibitem{Work} Workman, Myles. {\em Upper Semicontinuity of Index Plus Nullity for Minimal and CMC Hypersurfaces}. https://arxiv.org/abs/2312.09227.
 \end{thebibliography}
\end{document}